\documentclass[twoside,leqno]{article}

\usepackage[letterpaper]{geometry}

\usepackage{siamproceedings}

\usepackage[T1]{fontenc}
\usepackage{amsfonts}
\usepackage{graphicx}
\usepackage{epstopdf}
\usepackage{enumitem}
\usepackage{algorithmic}
\usepackage{xcolor}
\usepackage{dsfont}
\usepackage{mathrsfs}
\usepackage{wrapfig}
\usepackage{caption}
\usepackage{tikz}
\usepackage{tikz,tikz-cd}\usepackage{caption,subcaption}
\usetikzlibrary{arrows,decorations.pathmorphing,backgrounds,positioning,fit,petri}

\usetikzlibrary{decorations.pathreplacing}
\usepackage{subdepth}

\ifpdf
  \DeclareGraphicsExtensions{.eps,.pdf,.png,.jpg}
\else
  \DeclareGraphicsExtensions{.eps}
\fi

\newsiamremark{remark}{Remark}
\newsiamremark{assumption}{Assumption}

\usepackage{amsopn}

\newcommand{\bbE}{{\ensuremath{\mathbb E}} }

\newcommand{\bbR}{{\ensuremath{\mathbb R}} }

\newcommand{\bbV}{{\ensuremath{\mathbb V}} }

\newcommand{\cN}{{\ensuremath{\mathcal N}} }

\newcommand{\cU}{{\ensuremath{\mathcal U}} }

\newcommand{\cZ}{{\ensuremath{\mathcal Z}} }

\DeclareMathSymbol{\leqslant}{\mathalpha}{AMSa}{"36} 
\DeclareMathSymbol{\geqslant}{\mathalpha}{AMSa}{"3E} 
\DeclareMathSymbol{\eset}{\mathalpha}{AMSb}{"3F}     

\newcommand{\R}{\mathbb{R}}

\newcommand{\Z}{\mathbb{Z}}
\newcommand{\N}{\mathbb{N}}

\def\bs{\boldsymbol}

\newcommand{\PEfont}{\mathbf}

\newcommand{\p}{\ensuremath{\PEfont P}}

\DeclareMathOperator{\e}{\ensuremath{\PEfont E}}

\newcommand{\E}{\ensuremath{\PEfont  E}}
\renewcommand{\P}{\p}

\newcommand\bP{\ensuremath{\bs{\mathrm{P}}}}
\newcommand\bE{\ensuremath{\bs{\mathrm{E}}}}

\DeclareMathOperator{\bbcov}{\ensuremath{\mathbb{C}ov}}

\newcommand{\ind}{\mathds{1}}

\newcommand{\eps}{\varepsilon}
\renewcommand{\epsilon}{\varepsilon}
\renewcommand{\theta}{\vartheta}
\renewcommand{\rho}{\varrho}

\newcommand{\dd}{\mathrm{d}}

\newcommand{\cg}{\mathrm{cg}}

\begin{document}

\newcommand\relatedversion{}
\renewcommand\relatedversion{\thanks{The full version of the paper can be accessed at \protect\url{https://arxiv.org/abs/0000.00000}}} 

\title{\Large From disordered systems to the Critical 2D Stochastic Heat Flow}
    \author{Francesco Caravenna \thanks{University of Milano-Bicocca  (\email{francesco.caravenna@unimib.it}, \url{https://fcaraven.github.io}).}
    \and Rongfeng Sun \thanks{National University of Singapore  (\email{matsr@nus.edu.sg},
    \url{https://sites.google.com/view/rongfengsun/home}).}
    \and Nikos Zygouras\thanks{University of Warwick
  (\email{n.zygouras@warwick.ac.uk},  \url{https://warwick.ac.uk/fac/sci/maths/people/staff/zygouras/}).}}

\date{}

\maketitle







\begin{abstract} We review our joint work on the scaling limits of disordered systems, linking the notion of disorder relevance/irrelevance to that of sub/super-criticality of singular SPDEs. This line of research culminated in the construction of the Critical 2D Stochastic Heat Flow (SHF), a universal process which provides a non-trivial solution to the Stochastic Heat Equation in dimension 2, a critical singular SPDE that lies beyond the reach of existing solution theories. The SHF also offers a rare example of a non-Gaussian scaling limit for a disordered system at its phase transition point in the critical dimension.
\end{abstract}

\section{Introduction.}
The main question we want to address in this article is how to make sense of the
solution $u(t,x)$ of the {\em stochastic heat equation} (SHE)
\begin{equation}\label{eq:SHE}
	\partial_t u(t,x) =\tfrac{1}{2} \Delta u(t,x) + \beta \, \xi(t,x) \, u(t,x) \,,\qquad t>0, \ x\in\R^d,
\end{equation}
in the critical dimension $d=2$. This is a {\em singular stochastic partial differential equation} (SPDE) due to the ill-defined term
$\xi(t,x)u(t,x)$, where $\xi(t,x)$ is a space-time white noise and the solution $u(t,x)$ is expected to be a generalised function when $d\geq 2$.

The SHE is a fundamental object from various perspectives: (i) it is a proto-typical singular SPDEs;
(ii) it is a universal model for diffusive systems in disordered media \cite{BS95}; (iii) it is connected to random interface growth models such as the celebrated Kardar-Parisi-Zhang (KPZ) equation via the Cole-Hopf transformation $h=\log u$ \cite{KPZ86,KS91, Cor12, QS15}.

It is a classical result that \eqref{eq:SHE} is well-posed in dimension $d=1$.
However,  $d=2$ is critical and solution theories break down.
One can understand why $d=2$ is critical by the following scaling argument: if we scale space and time diffusively and define $\tilde u(t, x) := u(\epsilon^2 t, \epsilon x)$
for some $\eps>0$, then formally $\tilde u$ solves the equation
\begin{align}\label{eq:SHErg}
\partial_{ t} \tilde u =\frac{1}{2} \Delta \tilde u +
\beta\, \epsilon^{1-\frac{d}{2}} \,\tilde \xi \, \tilde u \,,
\qquad \tilde t>0, \ \tilde x\in\R^d,
\end{align}
where $\tilde \xi(t, x) := \epsilon^{1+\frac{d}{2}} \xi(\epsilon^2 t, \epsilon x)$ has the same distribution as $\xi$.
Since the ill-posedness of singular SPDEs is due to the loss of regularity on small scales, we zoom into smaller and smaller space-time scales
by sending $\epsilon\downarrow 0$ and note that, the strength of the noise vanishes in dimensions $d<2$, diverges in dimensions $d>2$, and stays constant in the critical
dimension $d=2$. If the equation contains non-linear terms, one can similarly consider the effect of scaling on the non-linearity.
In the language of singular SPDEs \cite{H13, H14, GIP15}, the SHE in dimensions $d<2$, $d=2$ and $d>2$ are called respectively {\em sub-critical}, {\em critical}, and {\em super-critical}. In recent years, there have been tremendous progress in developing solution theories for sub-critical singular
SPDEs \cite{H13, H14, GIP15, K16, GJ14, D22}. However, results for critical and super-critical singular SPDEs remain limited \cite{CSZ20, CSZ23a,  CES21, CET23a, CET23b, CGT24, CMT25, GP25, DG22, DG23}.

\smallskip

The main goal of this article to review recent progress that eventually led to the construction of the {\em Critical 2D Stochastic Heat Flow} (SHF) \cite{CSZ23a}, which provides a non-trivial solution to the 2D SHE. As we will discuss, the SHF also provides a rare example of a non-Gaussian scaling limit for a statistical physics model (the Directed Polymer Model) at the critical dimension and at its phase transition point.

We actually came to the study of the 2D SHE (and the associated 2D KPZ equation) from a different direction, namely, that of {\it disorder relevant systems} in statistical mechanics and their scaling limits. It turns out that the notions of sub-criticality, criticality, and super-criticality for singular SPDEs correspond to the notions of so-called disorder relevance, marginal relevance/irrelevance, and disorder irrelevance for  systems given by a ``pure'' model perturbed by disorder. To draw a first analogy with SPDEs, we may view the SHE \eqref{eq:SHE} as a
perturbation of the ``pure'' heat equation $\partial_t u = \frac{1}{2}\Delta u$
by the ``disorder'' term $\beta\xi u$.

Given a pure model, say defined on a lattice, a disorder perturbation is called relevant
if any amount of disorder alters the model's large scale behavior. In this case, it is natural to tune
the disorder strength down to~$0$ as the lattice spacing tends to~$0$, such that it leads to a continuum model with non-trivial dependence on disorder. In \cite{CSZ17a}, we developed a general theory to construct such continuum limits of disorder relevant models (see
Section~\ref{sec:disorder}). However, this theory does not apply to {\em marginally relevant} disordered systems, such as the 2D \emph{Directed Polymer Model} (DPM).

Interestingly, there is a close connection between the DPM and the SHE due to the fact that the partition functions of the DPM are precisely the solution of a \emph{regularised SHE} obtained by discretising space and time, which is an alternative to the regularisation via mollification of the noise (see Section~\ref{S:DRDP}). Our investigation of the scaling limit of the 2D DPM partition functions eventually led to the construction of the critical 2D SHF.

\smallskip

The starting point of our analysis was the identification of a \emph{phase transition} for the 2D DPM and for the 2D SHE in an \emph{intermediate disorder regime} \cite{CSZ17b} (see Section~\ref{S:phase}). Below the phase transition point, in the so-called \emph{sub-critical regime},
the partition functions of the 2D DPM and the solution of the mollified SHE
converge to the solution of the deterministic heat equation, and their fluctuations
converge to a Gaussian field that solves an \emph{additive} stochastic heat equation, thus belonging to the {\em Edwards-Wilkinson universality class} \cite{BS95}. The solution of the regularised 2D KPZ equation (which is connected to the regularised 2D SHE via the Cole-Hopf transformation $h=\log u$) was also shown in \cite{CSZ20, G20} to have the same
Gaussian fluctuations (see Section~\ref{sec:subcritical}, where further results in the sub-critical regime are reviewed).

At the phase transition point, more precisely in a critical window around it,
we proved in \cite{CSZ23a} that
a \emph{unique random scaling limit} emerges from  the 2D DPM partition functions,
named the {\em Critical 2D Stochastic Heat Flow} (SHF)
 (see Section~\ref{sec:construction}).
The SHF is a family of measure valued stochastic processes $\mathscr{Z}^\theta$
indexed by a renormalised disorder strength parameter $\theta \in \R$,
which can be interpreted as non-trivial solutions of the 2D SHE.
It was later shown by Tsai \cite{T24}
that the solution of the mollified 2D SHE also converges to the SHF.

\smallskip

Many features of the SHF have been established recently (see Section~\ref{sec:properties}),
including an axiomatic characterisation \cite{T24} and various properties concerning its invariance, singularity, regularity, asymptotic behavior, moments, connection to Gaussian Multiplicative Chaos (GMC) and black noise/noise sensitivity
\cite{CSZ23b, CM24, CSZ25, GT25, CD25, CCR25, CT25, N25b, BCT25, GN25, N25, C25a}.
We refer to the lecture notes \cite{CSZ24} for a more extended review
and open problems.

\smallskip

In the rest of this article, we will review the results sketched above
from an historical and personal perspective.
The paper is structured as follows.
\begin{itemize}
\item Section~\ref{S:DRDP} provides a high-level discussion of the connection between
disordered systems and SPDEs, and between the DPM and the SHE.
\item Section~\ref{sec:disorder} presents our general results on the scaling limits of disorder relevant
systems.
\item Section~\ref{S:phase} describes the phase transition of the 2D DPM on an intermediate disorder scale.

\item Section~\ref{sec:subcritical} discusses results in the sub-critical regime,
including 2D SHE and KPZ.
\item Section~\ref{sec:construction} reviews the construction of the SHF from \cite{CSZ23a}.
\item Section~\ref{sec:properties} presents recent progress on the SHF.
 \end{itemize}

\medskip
\section{Disorder relevance/irrelevance, scaling limits, directed polymer, and SHE.} \label{S:DRDP}

We review here some background, including the notions of disorder relevance/irrelevance for disordered systems, how they correspond to the notions of sub/super-criticality for singular SPDEs, how disorder relevance naturally leads to the construction of non-trivial continuum limits, and how such continuum limits for the directed polymer provide solutions to the SHE.

\subsection{Disordered systems and disorder relevance/irrelevance.}

We focus on disordered systems that are disorder perturbations of an underlying pure model, defined via Gibbs measures. Heuristically, {\em disorder
relevance} means that, arbitrarily small amount of disorder (or random environment) changes the large scale qualitative behaviour of the model (e.g., changes in the critical exponents) and leads to a different scaling limit, while {\em disorder irrelevance} means such qualitative changes arise only when the disorder strength is large enough. At the boundary between these two regimes, typically at a critical dimension, the effect of disorder perturbation is more subtle and could be either {\em marginally relevant} or {\em irrelevant}.

A key example is the Directed Polymer Model (DPM), which is a disorder perturbation of the simple random walk $S=(S_n)_{n\geq 0}$ on $\Z^d$ through
the family of Gibbs measures
\begin{align}\label{eq:DPM}
	\bP^{\omega}_{N, \beta}(S)
	:=\frac{1}{Z_{N, \beta}^{\omega}}
	e^{\sum_{n=1}^N \{\beta \omega(n,S_n) - \lambda(\beta)\}} \, \bP(S) \,,
\end{align}
where $\omega:=(\omega(n, x))_{n\in \N, x\in \Z}$ are i.i.d.\ disorder variables indexed by the time-space lattice $\N\times\Z^d$, $\beta\geq 0$ is the disorder strength, $\lambda(\beta):=\log \E[e^{\beta \omega(n,x)}]$, and the normalizing constant
\begin{align} \label{eq:paf}
	Z_{N, \beta}^{\omega} &= \bE\bigg[ e^{\sum_{n=1}^N \{\beta \omega(n,S_n) - \lambda(\beta)\}}\, \bigg| S_0=0\bigg]
\end{align}
is known as the partition function. As $\beta$ crosses a critical value $\beta_c\in [0, \infty)$, the DPM undergoes a phase transition from diffusive behaviour for the random walk, with path delocalisation, to expected super-diffusive behaviour and path localisation (see \cite{C17, Z24} for more background and recent progress). In particular, the DPM with $\beta=\infty$ is exactly the {\em last passage percolation} model, which belongs to the KPZ universality class \cite{QS15, Cor12, Cor16, Z22}. It turns out that $\beta_c=0$ in dimension $d=1$ and $2$,
and $\beta_c>0$ in $d\geq 3$. This implies that disorder is relevant for the DPM in $d=1$ and $2$, since any $\beta>0$ changes the large
scale behaviour of the model, whereas disorder is irrelevant in $d\geq 3$.  Dimension $d=2$ turns out to be the critical dimension, and hence disorder is called {\em marginally relevant} for the 2D DPM.

The connection between singular SPDEs and disordered systems is that, the notions of disorder relevance, marginal relevance/irrelevance, and disorder irrelevance for disordered systems correspond 
to the notions of sub-criticality, criticality, and super-criticality for singular SPDEs. Indeed, if we zoom out to larger and larger space-time scales by sending $\eps\uparrow \infty$ in the rescaled SHE \eqref{eq:SHErg}, then the dependence of the noise strength $\beta\, \epsilon^{1-\frac{d}{2}}$ on the dimension $d$ shows that, the noise term is a relevant perturbation of the deterministic
heat equation in $d<2$, a marginal perturbation in $d=2$, and an irrelevant perturbation in $d>2$.

\subsection{Scaling limits of disorder relevant models.} \label{S:DR}

In the physics literature, the Harris criterion \cite{H74} predicts whether disorder is relevant or irrelevant based on a relation between the spatial
dimension and the correlation length exponent of the pure model without disorder. In the mathematics literature, the question of disorder relevance/irrelevance has been investigated extensively for the disordered pinning model (see \cite{GT06, AZ09, GLT10} and the references in
\cite{Gi11}), and also for the random field Ising model (RFIM) \cite{BK88, AW90, DX21, AHP20, DZ24} and the DPM \cite{L25b}.

Inspired by \cite{AKQ14a, AKQ14b} on the scaling limit of the DPM in dimension $1+1$, we proposed in \cite{CSZ17a} a new perspective
on disorder relevance/irrelevance. The heuristic is that, if the disorder perturbation of a pure model is relevant, then as we send the lattice spacing $\delta\downarrow 0$ (equivalent to zooming out to larger and larger space-time scales), the effect of disorder will diverge. It is then natural to compensate this by tuning the disorder strength down to $0$ at a suitable rate depending on $\delta$, such that in the limit we obtain a continuum model with non-trivial dependence on disorder. Note that this is only possible if disorder is relevant, and hence we can use the existence of a non-trivial
disordered continuum limit as a criterion for disorder relevance. This also motivates us to construct non-trivial scaling limits of disorder relevant models, which has been successfully carried out for the disordered pinning model \cite{CSZ16} and the random field perturbation of the critical 2D Ising model \cite{CSZ17a, BS22}.

We will review this theory in Section \ref{sec:disorder}, which is based on polynomial chaos expansions for the model's partition functions and their convergence to Wiener-It\^o chaos expansions. However, this theory is not applicable to marginally relevant systems such as the 2D DPM.

\subsection{Directed Polymer Model and Stochastic Heat Equation.}
\label{sec:SHE-DP}

We now explain the precise  connection between the DPM and the SHE, which allows us to interpret the continuum limit of the 2D DPM as a solution of the 2D SHE.

To make sense of singular SPDEs such as the SHE \eqref{eq:SHE}, which are ill-posed due to the irregularity of the noise $\xi$ and the solution
$u$ on small spatial scales, the standard procedure is to
first perform a regularisation on the spatial scale $\eps$ (also known as ultraviolet cutoff) and then take suitable limits
as $\eps\downarrow 0$. If the solution $u^\eps$ of the regularised equation admits a non-trivial limit after suitable tuning of parameters, centering
and scaling, then the limit can be interpreted as
a solution of the singular SPDE.

There are different ways of performing the ultraviolet cutoff, which are expected to lead to the same limits by universality. These include:
\begin{itemize}
\item[(a)] Mollify the space-time white noise $\xi$ on the spatial scale $\eps$;

\item[(b)] Discretize space (or space-time) on the spatial scale $\eps$;

\item[(c)] Truncate the Fourier transform at frequency $1/\eps$.
\end{itemize}
In approach (a), we replace $\xi$ by its spatial mollification $\xi^\eps := j^\epsilon * \xi$, where $j^\eps(x):=\epsilon^{-2}j(\frac{x}{\eps})$ for
some smooth probability density $j\in C_c^\infty(\R^2)$ with compact support. This leads to the mollified 2D SHE
\begin{equation} \label{eq:mollSHE}
\partial_t u^\eps = \frac{1}{2}\Delta u^\eps + \beta_\eps u^\eps \xi^\eps, \qquad u^\eps(0, \cdot)\equiv 1,
\end{equation}
which admits a classical It\^o solution and a Feynman-Kac representation (see \cite[Section 3]{BC95}):
\begin{equation} \label{eq:FK1}
u^\eps(t,x) = \bE\Big[e^{\int_0^t (\beta_\eps \xi^\eps(t-s, B_s) - \lambda_\eps){\rm d}s} \Big| B_0=x \Big],
\end{equation}
where $\bE[\cdot]$ is w.r.t.\ a standard Brownian motion $B$ in $\R^2$, and $\lambda_\eps:= \frac{\beta_\eps^2}{2} \bbV{\rm ar}(\xi^\eps(s, x))=\frac{\beta_\eps^2}{2 \eps^2}\Vert j\Vert_2^2$ ensures that $\bbE[u^\eps(t,x)]=1$. Apart from a time reversal, this is the continuum analogue of the DPM partition function in \eqref{eq:paf}, which leads to approach (b).

If we follow (b) and discretize time and space, then the Brownian motion $B$ in \eqref{eq:FK1} is replaced by a simple random walk, the mollified space-time white noise $\xi^\eps$ is replaced by i.i.d.\ disorder on the diffusive rescaled time-space lattice $\eps^2 \N \times \eps \Z^2$, and we have
the following correspondence between $u^\eps$ and the time-space rescaled DPM partition functions
\begin{equation}\label{eq:ueps-ZN}
	u^\eps(1-t, x) \longleftrightarrow Z_{N, \beta}^{\omega}([Nt], [\sqrt{N} x]),
	\qquad t\in [0, 1], x\in \R^2,
\end{equation}
where $\eps=\frac{1}{\sqrt N}$, brackets $[\cdot]$ denotes integer parts, and $Z_{N, \beta}^{\omega}(m, z)$ is the point-to-plane partition function from the time-space point $(m, z)\in \N\times\Z^2$ to the plane $\{N\}\times \Z^2$:
\begin{align} \label{eq:paf2}
	Z_{N, \beta}^{\omega}(m, z) &= \bE\bigg[ e^{\sum_{n=m+1}^N \{\beta \omega(n,S_n) - \lambda(\beta)\}}\, \bigg| S_m=z\bigg].
\end{align}

The problem of defining a non-trivial solution of the 2D SHE now becomes equivalent to
finding non-trivial scaling limits for the family of
rescaled DPM partition functions $(Z_{N, \beta}^{\omega}([Nt], [\sqrt{N} x]))_{t\in [0, 1], x\in \R^2}$. The construction of the limit, the {\em Critical 2D Stochastic Heat Flow} \cite{CSZ23a} (see Sections~\ref{sec:construction} and~\ref{sec:properties}), builds on a deeper understanding of the 2D DPM, and in particular, starts with the discovery of a phase transition for the 2D DPM in an intermediate disorder regime \cite{CSZ17b}. We will review these results in Section \ref{S:phase}. In the next Section~\ref{sec:disorder}, we will first review the theory on the scaling limits of disorder relevant models.

\medskip
\section{Scaling limits of disorder relevant systems.}\label{sec:disorder}
As explained in Section \ref{S:DR}, when the disorder perturbation of a pure model is relevant for the model's large scale behaviour, then it should be possible to tune disorder strength down to $0$ as a function of the lattice spacing, such that we obtain a continuum limit with non-trivial dependence on disorder. For disordered systems defined via Gibbs measures, the existence of such continuum limits should first manifest itself through the convergence of the partition functions. We now review the theory developed in \cite{CSZ17a} for the convergence of partition functions of disorder relevant systems, which was inspired by a result of this nature for the 1D DPM \cite{AKQ14a}.

The class of disordered systems we focus on are random field perturbations of binary-valued spin systems. Let $\Omega \subset \R^d$ be an open and simply connected domain. Given lattice spacing $\delta>0$, let $\Omega_{\delta}:=\Omega \cap (\delta \mathbb{Z})^{d}$ denote the lattice approximation of $\Omega$. More generally, the lattice spacing could be different in different directions and we could define $\Omega_{\delta}:=\Omega \cap (\delta^{a_1}\Z \times \delta^{a_2}\Z\times \cdots \times\delta^{a_d}\Z)$ for some $a_1, \ldots, a_d>0$.
The pure (or reference) model will be a spin system $\sigma=(\sigma_x)_{x\in \Omega_\delta}$ with a spin $\sigma_x\in \{0, 1\}$ assigned to each $x\in \Omega_\delta$. Let $\p_{\Omega_{\delta}}^{\text {ref }}$ (with expectation $\e_{\Omega_{\delta}}^{\text {ref }}$) denote the law of $\sigma$. We will assume that $\p_{\Omega_{\delta}}^{\text {ref }}$ has a continuum limit as $\delta\downarrow 0$, which is expected to hold for equilibrium spin systems at the critical point of a continuous phase transition.

The disorder we consider will be i.i.d.\ random fields $\omega=(\omega_x)_{x\in \Omega_\delta}$ with $\bbE[\omega_x]=0$, $\bbE[\omega_x^2]=1$,
and finite moment generating function $\lambda(\beta):= \log \bbE[e^{\beta \omega_x}]$ for all $\beta$ in a neighbour around $0$.

Given disorder $\omega$, field bias $h\in \R$, and disorder strength $\beta>0$, we can define the random field perturbation of
$\p_{\Omega_{\delta}}^{\text {ref }}$ via the Gibbs measure
\begin{equation}\label{eq:Gibbs}
\p_{\Omega_{\delta} ; \beta, h}^{\omega}(\mathrm{d} \sigma):=\frac{e^{\sum_{x \in \Omega_{\delta}} (\beta \omega_{x} +h) \sigma_{x}}}{Z_{\Omega_{\delta} ; \beta, h}^{\omega}} \bP_{\Omega_{\delta}}^{\mathrm{ref}}(\mathrm{d} \sigma),
\end{equation}
where the partition function is defined by
\begin{equation}\label{eq:Zom}
Z_{\Omega_{\delta} ; \beta, h}^{\omega}:=\E_{\Omega_{\delta}}^{\mathrm{ref}}\left[e^{\sum_{x \in \Omega_{\delta}} (\beta \omega_{x} +h)\sigma_{x}}\right] .
\end{equation}
If such a random field perturbation of the pure model $\p_{\Omega_{\delta}}^{\text {ref }}$ is relevant, then it should be possible to find an {\em intermediate disorder regime} with  $\beta=\beta_\delta \downarrow 0$ and $h=h_\delta \downarrow 0$, such that after proper centering and scaling,
the partition function $Z_{\Omega_{\delta} ; \beta_\delta, h_\delta}^{\omega}$ admits a non-trivial distributional limit.

Examples of disordered systems that fit within this framework include the following:

\begin{itemize}
\item {\bf 1D Directed Polymer Model.} For the directed polymer with partition function as in \eqref{eq:paf}, the spin field on $\N\times \Z$ is defined from the random walk path $S$ by $\sigma_{(n, x)}:= 1_{\{ S_n=x\}}$. Under diffusive rescaling of the time-space lattice, we see that $\Omega = (0,1) \times \R$ and $\Omega_\delta :=  \Omega\cap (\delta \Z \times \delta^{1/2}\Z)$, with $\delta = N^{-1}$.

\item {\bf Critical 2D Ising Model.} The reference measure is that of the critical 2D Ising model with $+$ boundary condition, where $\Omega\subset \R^2$ is an open, bounded, simply connected domain, $\Omega_\delta =\Omega \cap \delta \Z^2$, $\sigma_x\in \{\pm 1\}$, and the reference measure is defined by
\begin{align*}
\P^{\rm ref}_{\Omega_\delta}(\sigma):=
\frac{1}{Z^+_{\Omega_\delta; \beta_c}}e^{\sum_{x\sim y \in \Omega_\delta \cup \partial \Omega_\delta} \beta_c \sigma_x \sigma_y},
\end{align*}
where the sum is over all edges between vertices in $\Omega_\delta\cup \partial \Omega_\delta$, $\beta_c=\log (1+\sqrt{2})/2$ is the critical inverse temperature of the 2D Ising model, $\partial\Omega_\delta$ is the outer
boundary of  $\Omega_\delta$, and $\sigma_y\equiv 1$ for all $y\in \partial \Omega_\delta$. A linear change of variable $(\sigma_x+1)/2\in \{0, 1\}$ recasts the reference measure into the framework above.
\end{itemize}

We now recall from \cite{CSZ17a} the general convergence criteria for the disordered partition functions $Z_{\Omega_{\delta} ; \beta_\delta, h_\delta}^{\omega}$ to have non-trivial continuum limits (an extension to non-binary spins can be found in \cite{LMS24}).

For $k\in \N$ and $x_1, \ldots, x_k \in \Omega$, define the \emph{$k$-point correlation function} of the reference measure $\p_{\Omega_{\delta}}^{\text {ref }}$ by
\begin{equation} \label{eq:corrfun}
	\psi_{\Omega_\delta}^{(k)}(x_1, \ldots, x_k) :=
	\begin{cases}
	\E_{\Omega_\delta}^{\rm ref} \big[\sigma_{x_1^{(\delta)}} \, \sigma_{x_2^{(\delta)}}
	\cdots \sigma_{x_k^{(\delta)}} \big]
	& \text{if } x_i^{(\delta)} \ne x_j^{(\delta)} \text{ for all } i\ne j , \\
	0 & \text{otherwise},
	\end{cases}
\end{equation}
where $x^{(\delta)}$ denotes the point in $\Omega_\delta$ closest to $x$.

\begin{assumption}\label{ass:main}
There exists a correlation exponent $\gamma \in [0,\infty)$ such that, for all $k\in\N$,
there exists a symmetric function $\boldsymbol{\psi}_\Omega^{(k)}: \Omega^k \to \R$ and
\begin{equation} \label{eq:polyconv}
	(\delta^{-\gamma})^k \, \psi_{\Omega_\delta}^{(k)} (x_1, \ldots, x_k)
	\,\xrightarrow[\,\delta \downarrow 0\,]{}\,
	\boldsymbol{ \psi}_{\Omega}^{(k)} (x_1, \ldots, x_k) \qquad
       \text{in } \,\,L^2(\Omega^k) \,.
\end{equation}
Furthermore, for some $\eps>0$,
\begin{equation} \label{eq:L2sum}
	\limsup_{\ell\to\infty} \limsup_{\delta\downarrow0}\sum_{k>\ell}
	\frac{(1+\eps)^k}{k!} \, \big\| \psi_{\Omega_\delta}^{(k)} \big\|_{L^2(\Omega^k)}^2 = 0.
\end{equation}
\end{assumption}

We can now state the convergence result from \cite{CSZ17a}.
 \begin{theorem}\label{th:mainhere}
Let $Z_{\Omega_{\delta} ; \beta_\delta, h_\delta}^{\omega}$ be the disordered partition function defined as in \eqref{eq:Zom}.
Suppose the reference measure $\p_{\Omega_{\delta}}^{\text {ref }}$ satisfies Assumption~\ref{ass:main} for some exponent
$\gamma<d/2$. Then choosing
\begin{equation}\label{eq:bh}
\beta_\delta =\hat\beta \delta^{\frac{d}{2}-\gamma}, \qquad h_\delta= - \log \bbE[e^{\beta_\delta \omega_x}] +  \hat h \delta^{d-\gamma}
\qquad \mbox{for some } \hat\beta>0 \mbox{ and } \hat h\in \R \,,
\end{equation}
we have
\begin{equation*}
	Z_{\Omega_\delta; \beta_\delta, h_\delta}^\omega
	\ \xrightarrow[\,\delta \downarrow 0\,]{(d)} \
	\mathscr{Z}^W_{\Omega; \hat\beta, \hat h},
\end{equation*}
where $W$ is a white noise on $\Omega$, and $\mathscr{Z}^W_{\Omega; \hat\beta, \hat h}$ admits the following {\it Wiener-It\^o chaos expansion}
\begin{align}\label{eq:ZW}
\mathscr{Z}^W_{\Omega; \hat\beta} := 1 +
	\sum_{k=1}^{\infty} \;
	\frac{1}{k!} \,
	\idotsint_{\Omega^k}
	\boldsymbol{\psi}^{(k)}_{\Omega}(x_1, \ldots, x_k)
	\, \prod_{i=1}^k \, (\hat \beta W(\text{\rm d} x_i) +\hat h {\rm d}x_i) \,.
\end{align}
\end{theorem}
\begin{remark}
Disorder relevance enters through the condition $\gamma<d/2$ in Theorem \ref{th:mainhere}, because otherwise $\beta_\delta=\hat\beta \delta^{d/2-\gamma}$ will diverge and it will not be possible to obtain a non-trivial limit for the partition function by sending $\beta\downarrow 0$. This condition on $\gamma$
is consistent with the Harris criterion (see \cite[Section 1.3]{CSZ17a} for more details). When $\gamma=d/2$, which is the marginal case, this approach fails because even though the limiting correlation function $\boldsymbol{ \psi}_{\Omega}^{(k)}$ in \eqref{eq:polyconv} may still exist pointwise, it just fails to be square integrable, and hence the stochastic integrals in \eqref{eq:ZW} are undefined. This will be the case for the partition function of the 2D DPM.
\end{remark}

\begin{remark}
The proof of Theorem \ref{th:mainhere} starts with a discrete analogue of the Wiener-It\^o chaos expansion for the partition function $Z_{\Omega_{\delta} ; \beta_\delta, h_\delta}^{\omega}$ with respect to the disorder, known as {\em polynomial chaos expansion}. More precisely,
in \eqref{eq:Zom}, we can use the fact $\sigma_x\in \{0, 1\}$ to rewrite
$e^{(\beta \omega_x + h)\sigma_x} = 1 +\eta_x \sigma_x$ with
$\eta_x:=  e^{\beta \omega_x + h}-1$.
Expanding the product over $x\in \Omega_\delta$ in \eqref{eq:Zom} then leads to the following polynomial chaos expansion in the i.i.d.\ variables $(\eta_x)_{x\in \Omega_\delta}$:
 \begin{equation} \label{eq:Zstart2}
	Z_{\Omega_{\delta} ; \beta_\delta, h_\delta}^{\omega} =
	1 + \sum_{k=1}^{|\Omega_\delta|} \;
	\frac{1}{k!} \,
	\sum_{(x_1, x_2, \ldots, x_k) \in (\Omega_\delta)^k}
	\psi^{(k)}_{\Omega_\delta}(x_1, \ldots, x_k)
	\, \prod_{i=1}^k \eta_{x_i} .
\end{equation}
As $\beta=\beta_\delta, h=h_\delta\to 0$ as $\delta\downarrow 0$, the influence of each $\eta_x$ on $Z_{\Omega_{\delta} ; \beta_\delta, h_\delta}^{\omega}$ becomes asymptotically negligible, which allows us to apply a Lindeberg principle \cite{Cha06, MOO10} to replace $(\eta_x)_{x\in \Omega_\delta}$ by a family of i.i.d.\ Gaussian random variables $(\xi_x)_{x\in \Omega_\delta}$ with the same mean and variance. Furthermore,
a polynomial chaos expansion in $(\xi_x)_{x\in \Omega_\delta}$ can be written as a Wiener-It\^o chaos expansion w.r.t.\ a white noise $W$ on $\Omega$. The choice of $\beta_\delta$ and $h_\delta$ in \eqref{eq:bh} and Assumption \ref{ass:main} then ensure that this sequence of Wiener-It\^o chaos expansions will converge.
\end{remark}

\begin{remark}
It is possible to extend the convergence of disordered partition functions in Theorem \ref{th:mainhere} to construct non-trivial continuum limits of the random Gibbs measure $\p_{\Omega_{\delta} ; \beta, h}^{\omega}(\mathrm{d} \sigma)$, where the disorder $\omega$ on $\Omega_\delta$ converges to a white noise $W$ on $\Omega$. This has been carried out for the DPM in dimension $1+1$ \cite{AKQ14b}, the disordered pinning model \cite{CSZ16}, which can be regarded as a DPM in dimension $1+0$, with the random walk $S$ on $\Z^d$ replaced by a renewal process. It has also been carried out for the random field perturbation of the critical 2D Ising model \cite{BS22}. The key observations behind these constructions is that, for the directed models, the measure $\p_{\Omega_{\delta} ; \beta, h}^{\omega}(\mathrm{d} \sigma)$ is determined by the joint law of the family of point-to-point partition functions, while for the random field Ising model, the law of the spin field is determined by its Fourier transform, which can be written in terms of the Ising partition function with an imaginary external field.
\end{remark}

\section{Phase transition in the 2D DPM.}
\label{S:phase}

We now address the question of finding non-trivial limits for the partition functions of the 2D DPM, which are solutions of a regularised version of the 2D SHE as discussed in Section \ref{sec:SHE-DP}. The starting point is the discovery that even though the 2D DPM has critical value $\beta_c=0$, there is still a phase transition in an intermediate disorder regime \cite{CSZ17b}.

Recall from \eqref{eq:DPM}-\eqref{eq:paf} the directed polymer measure $\bP^{\beta,\, \omega}_{N}(S)$, with partition function
$$
	Z_{N, \beta}^{\omega} = \bE\bigg[ e^{\sum_{n=1}^N \{\beta \omega(n,S_n) - \lambda(\beta)\}}\, \bigg| S_0=0\bigg].
$$
Since $\lambda(\beta)=\log \bbE[e^{\beta\omega(n, x)}]$, it is easily seen that $(Z_{N, \beta}^{\omega})_{N\in\N}$ is a non-negative martingale. The critical point $\beta_c(d)$ is defined as the boundary between a {\em weak disorder} phase where $Z_{N, \beta}^{\omega}\to Z_\infty^{\beta, \omega}>0$ almost surely for $\beta < \beta_c$, and a {\em strong disorder} phase where  $Z_{N, \beta}^{\omega}\to 0$ almost surely for $\beta > \beta_c$
(see \cite{C17} and the references therein).
It is known that $\beta_c=0$ in $d=1, 2$, and $\beta_c\in (0, \infty)$ in  $d\geq 3$. In a recent breakthrough \cite{JL24, JL25}, it was shown for $\beta = \beta_c$ in $d\geq 3$, the DPM also belongs to the weak disorder phase. In the weak disorder phase, it is known that the polymer is diffusive and satisfies an invariance principle \cite{CY06, L25}, and hence the polymer is delocalised in space.  On the other hand, in the strong disorder phase, the partition function $Z_{N, \beta}^{\omega} \to 0$ exponentially fast in $N$ \cite{JL24, JL25}, which is
known to imply a form of path localisation $\omega$ \cite{CH02, V07, Cha19, BC20a, BC20b}.

In the critical dimension $d=2$, it was discovered in \cite{CSZ17a} that there is still a phase transition in an intermediate disorder regime $\beta_N=\hat\beta /\sqrt{R_N}$ for some diverging sequence $R_N$, such that for $\hat\beta< \hat\beta_c=1$, the partition function converges in law to a positive limit, while for $\hat\beta\geq \hat\beta_c=1$, the partition function converges in law to $0$. Such a dichotomy is precisely what characterises the phase transition of the DPM in $d\geq 3$.

We now give the precise statement of the result. Let
\begin{align}\label{2interm}
\beta_N:=\frac{\hat\beta}{\sqrt{R_N}}
 \qquad \text{with} \qquad R_N:= \sum_{n=1}^N \bP(S_n=S'_n) = \frac{1+o(1)}{\pi} \log N,
 \end{align}
 where $S$ and $S'$ are two i.i.d.\ simple random walks on $\Z^d$.

\begin{theorem}[\cite{CSZ17b}]\label{2dsubcritical}
Let $Z_{N,\beta_N}^{\omega}$ be the partition function of the 2D DPM with $\beta_N$ chosen as in \eqref{2interm}.
 We have
 \begin{equation} \label{conv-subZ}
	Z_{N,\beta_N}^\omega \xrightarrow[\,N\to\infty\,]{dist}
	\begin{cases}
	\exp\big( \sigma_{\hat\beta} \, \cN-
	\frac{1}{2} \sigma_{\hat\beta}^2
	\,\big)
	& \text{if } \hat \beta < 1 \\
	0 & \text{if } \hat\beta \ge 1
	\end{cases} \,,
\end{equation}
where $\cN$ is a standard normal variable and $ \sigma_{\hat\beta}^2  = \log (1-\hat\beta^2)^{-1}$.
\end{theorem}

Since the mollified SHE \eqref{eq:mollSHE} and the DPM partition functions give alternative regularisations of the 2D SHE (see Section \ref{sec:SHE-DP}), it should come as no surprise that an analogue of Theorem~\ref{2dsubcritical} also holds for the mollified 2D SHE.

\begin{theorem}[\cite{CSZ17b}]\label{2dsubcriticalSHE}
Let $u^\eps(t, x)$ be the solution of the mollified 2D SHE  \eqref{eq:mollSHE} with disorder strength
\begin{equation}\label{eq:betaeps}
\beta_\epsilon :=(\hat\beta +o(1))\sqrt{\frac{2\pi}{\log\frac{1}{\epsilon}}}.
\end{equation}
Then for every $(t, x)\in (0, \infty)\times \R^2$, we have
 \begin{equation} \label{conv-subZ-SHE}
	u^\eps(t, x) \xrightarrow[\,N\to\infty\,]{dist}
	\begin{cases}
	\exp\big( \sigma_{\hat\beta} \, \cN-
	\frac{1}{2} \sigma_{\hat\beta}^2
	\,\big)
	& \text{if } \hat \beta < 1 \\
	0 & \text{if } \hat\beta \ge 1
	\end{cases} \,.
\end{equation}
\end{theorem}

\begin{remark}\label{rem:further}
In \cite{CSZ17b}, it was also shown that $u^\eps(t, x)$ at distinct time-space points converge to independent limits. Furthermore, the same phase transition and log-normal limit were established for two other marginally relevant polymer models. One is a marginally relevant disordered pinning model, and the other is the 1D DPM defined from Cauchy random walks. The proof relies on an intrinsic exponential separation of space-time scales common to all these models.
\end{remark}

\begin{remark}
Recently, Theorem \ref{2dsubcriticalSHE} has been extended in \cite{DG22, DG23, DG24} to 2D semi-linear SHE's of the form
$$
 \partial_t u^\epsilon=\frac{1}{2} \Delta u^\epsilon+
  (\log \tfrac{1}{\epsilon})^{-\frac{1}{2}} \sigma(u^\epsilon) \,\xi^\epsilon,
$$
where $\sigma$ is a Lipschitz function. It was shown that under suitable assumptions on $\sigma$, $u^\eps(t, x)$ converges in distribution
to the terminal value of an associated forward-backward SDE, which includes Theorem  \ref{2dsubcriticalSHE}  as a special case.

\end{remark}

\begin{remark}
A particular feature of 
dimension 2, 
in contrast to higher dimensions $d\geq 3$,
is that the critical point in the former case is 
the point at which the
second moment $\bbE\big[ \big(Z_{N,\beta_N}^\omega)^2 \big]$ blows up as $N\to\infty$.
A similar phenomenon holds for the 2D SHE and other marginally relevant polymer models
\cite{CSZ17b}.
Furthermore, we have that below the critical point in dimension $d=2$, {\it all} moments remain finite
\cite{LZ23, CZ23, LZ24, CZ24}, while at the critical point, the $h$-moments, for any fixed $h\geq 2$,
blow up as $(\log N)^{\frac{h(h-1)}{2}+o(1)}$ as demonstrated in a continuous setting in \cite{LiuZ24}.
Contrary to this sharp transition, the directed polymer in $d\geq 3$ admits a gradual loss of moment
as the critical point is approached from below \cite{J22}.
\end{remark}

\medskip
\section{2D SHE and KPZ in the subcritical regime.}
\label{sec:subcritical}

From the point of view of singular SPDE's, we are interested in the solution of the regularised 2D SHE (and 2D KPZ) in the random field limit (regarded as random generalised functions). It turns out that the critical point identified in Theorems \ref{2dsubcritical} and \ref{2dsubcriticalSHE} for the one-point distribution is also the critical point for the random field limit. Below the critical point, the solution of the regularised 2D SHE converges to the solution of the additive stochastic heat equation, also known as the {\em Edwards-Wilkinson equation}. Surprisingly, the solution of the regularised 2D KPZ below the critical point converges to the same limit.

\subsection{Gaussian limit for the subcritical 2D SHE}

The convergence of the regularised solution of the 2D SHE in the subcritical regime was proved in \cite{CSZ17b} for both the DPM partition functions and the solution of the mollified SHE. To avoid repetitions, we will only state the result for the mollified SHE.

\begin{theorem}[Edwards-Wilkinson fluctuations \cite{CSZ17b}]\label{T:EWDSHE}
Let $u^\eps$ be the solution of the mollified 2D SHE  \eqref{eq:mollSHE} with disorder strength $\beta_\eps=(\hat\beta+o(1)) \sqrt{\frac{2\pi}{\log \frac{1}{\epsilon}}}$ for some $\hat\beta\in (0,1)$. Then for any test function $\phi\in C_c(\R^2)$, we have
\begin{equation} \label{eq:EW-sub}
\frac{1}{\beta_\epsilon} \int_{\R^2} \phi(x) \big(u^\eps(t,x)-1\big) {\rm d}x
\xrightarrow[\,\eps\to 0\,]{d} \int_{\R^2} \phi(x) v(t, x) {\rm d}x,
\end{equation}
where $v(t, x)$ is the solution to the Edwards-Wilkinson equation
\begin{equation}\label{linearSHE}
\begin{aligned}
	\partial_t v(t,x) & =\frac{1}{2} \Delta v(t,x) + \sqrt{\frac{1}{1-\hat\beta^2}}\, \tilde\xi(t,x), \qquad v(0, x)\equiv 0\,.
\end{aligned}
\end{equation}
\end{theorem}
\begin{remark}
The space-time white noise $\tilde \xi$ in \eqref{linearSHE} is a mixture of the driving white noise $\xi$ in the mollified SHE  \eqref{eq:mollSHE}
and an independent white noise. See \cite{CSZ24} for more details. We also note that the noise coefficient in \eqref{linearSHE} diverges as
$\hat\beta\uparrow \hat\beta_c=1$.
\end{remark}

\begin{remark}\label{rem:LLN}
It follows from \eqref{eq:EW-sub} that, in the sub-critical regime $\hat\beta < 1$,
the mollified  solution $u^\eps$ satisfies a \emph{law of large numbers}, namely
$\int_{\R^2} \phi(x) \, u^\epsilon(t,x) \, \dd x \to \int_{\R^2}
\phi(x) \, \dd x$ in probability as $\epsilon \to 0$, for any $\phi\in C_c(\R^2)$.
\end{remark}

\subsection{Gaussian limit for the subcritical 2D KPZ}
The KPZ equation is formally given by
\begin{equation}\label{eq:KPZ}\tag{KPZ}
    \partial_t h(t, x) = \frac{1}{2} \Delta h(t, x) + \frac{1}{2} |\nabla h(t,x)|^2 +
    \beta\,\xi(t, x) \,,\qquad t>0, \ x\in\R^d.
\end{equation}
It is a model for random interface growth and has been studied extensively in $d=1$ as a canonical example in the KPZ universality class \cite{KPZ86,BS95, Cor12,QS15, Cor16, Z22}. The solution $h$ of the KPZ equation is related to the solution $u$ of the SHE by the Cole-Hopf transformation $h=\log u$. It is therefore not surprising that $d=2$ is also critical for the KPZ equation. We consider the solution $h^\eps$ of the mollified 2D KPZ, where the space-time white noise $\xi$ is replaced by its mollification $\xi^\eps$ as in the mollified 2D SHE \eqref{eq:mollSHE}. It can be seen that
$h^\eps(t, x):= \log u^\eps(t, x)$ solves the following mollified 2D KPZ equation
\begin{align} \label{eq:2DKPZ}
\partial_t h^\epsilon = \frac{1}{2}\Delta h^\epsilon + \frac{1}{2} |\nabla h^\epsilon|^2 +\beta_\epsilon \xi^\epsilon -
C_\epsilon \qquad \text{with} \quad  C_\epsilon:=\beta_\epsilon^2 \epsilon^{-2} \|j\|_2^2.
\end{align}

A surprising result is that, $h^\eps =\log u^\eps$ has the same random field limit as $u^\eps$ in the sub-critical regime, even though
$u^\eps(t, x)$ has a log-normal distribution that is not close to $1$.

\begin{theorem}[Edwards-Wilkinson fluctuations \cite{CSZ20, G20}]\ \label{thm:KPZ}
Let $h^\epsilon$ be the solution of the mollified 2D KPZ equation \eqref{eq:2DKPZ} with disorder strength $\beta_\eps=(\hat\beta+o(1)) \sqrt{\frac{2\pi}{\log \frac{1}{\epsilon}}}$ for some $\hat\beta\in (0,1)$. Then for any test function $\phi\in C_c(\R^2)$, we have
\begin{align}\label{EW-KPZ}
\frac{1}{\beta_\epsilon} \int_{\R^2} \phi(x) \Big(h^\eps(t,x)-\E[h^\eps(t,x)] \Big) {\rm d}x
\xrightarrow[\,\eps\to 0\,]{d} \int_{\R^2} \phi(x) v(t, x) {\rm d}x,
\end{align}
where $v(t, x)$ is the solution of the Edwards-Wilkinson equation \eqref{linearSHE}.
\end{theorem}
Theorem \ref{thm:KPZ} was proved in \cite{CSZ20} based on chaos expansion methods, and independently in \cite{G20} for $\hat\beta$
small using Malliavin calculus.

\begin{remark}
In a series of works \cite{CES21, CET23a, CET23b, CGT24}, it was shown that the {\it anisotropic} 2D KPZ equation
\begin{align}\label{AKPZ}
\partial_t h^\epsilon = \frac{1}{2}\Delta h^\epsilon + \frac{1}{2} \big((\partial_x h^\epsilon)^2 -(\partial_y h^\epsilon)^2 \big) +\frac{\hat\beta}{ (\log \eps^{-1})^{1/2}}  \xi^\epsilon \end{align}
also sees Edwards-Wilkinson fluctuations in the limit $\eps\downarrow 0$. However, the noise coefficient in the limiting EW equation is finite for all values of $\hat\beta$, while the diffusion coefficient also depends on $\hat\beta$.
The methods employed for the anisotropic 2D KPZ are very different from the (isotropic) 2D KPZ, and it will be interesting to understand the transition  between the two cases.
\end{remark}

\begin{remark}
An interesting work \cite{CNZ25} has initiated the study of the asymptotic maxima of $h^\epsilon$ in the sub-critical regime,
making connections but also drawing contrasts
to the theory of extrema of log-correlated fields. We refer to \cite{CNZ25}
for details.
\end{remark}

\medskip
\section{The Critical 2D Stochastic Heat Flow.}
\label{sec:construction}

We now go back to our original question: how can we construct a non-trivial solution to the SHE
\eqref{eq:SHE} in dimension $d=2$?
As we discussed in Section~\ref{sec:SHE-DP},
discretising space-time turns this question into the problem of \emph{finding
a non-trivial scaling limit for the rescaled partition functions
$(Z_{N, \beta}^{\omega}([Nt], [\sqrt{N} x]))_{t\in [0, 1], x\in \R^2}$} as $N\to\infty$. Solving this problem
leads to the construction of the \emph{critical 2D Stochastic Heat Flow (SHF)} \cite{CSZ23a},
which we describe in the present section.

\subsection{Construction of the SHF.}

As shown in Theorem \ref{T:EWDSHE} and Remark \ref{rem:LLN}, in the sub-critical regime $\hat\beta<\hat\beta_c=1$ on the
intermediate disorder scale defined in \eqref{2interm} and \eqref{eq:betaeps}, the  2D polymer partition functions $Z_{N, \beta_N}^{\omega}([Nt], [\sqrt{N} x])$ and  the solution $u^\eps(t, x)$ of the mollified SHE satisfy a law of large numbers, with Edwards-Wilkinson
fluctuations that blow up as $\hat\beta\uparrow 1$. It is then natural to investigate what happens at $\hat\beta=1$, whether the random field
$Z_{N, \beta_N}^{\omega}([Nt], [\sqrt{N} \cdot ])$ admits a non-trivial limit after suitable centering and scaling.

The first clue came from moment calculations for
$$
\cZ^{\beta_N}_{N;\, t}(\varphi) := \int_{\R^2} \varphi(x) \, Z_{N, \beta_N}^{\omega}([Nt], [\sqrt{N} x]) \, \dd x, \qquad \varphi \in C_c(\R^2).
$$
Clearly, the first moment satisfies $\bbE[\cZ^{\beta_N}_{N;\, t}(\varphi)] \to \int \varphi(x) \dd x$. It turns out that the second moment $\bbE[\cZ^{\beta_N}_{N;\, t}(\varphi)^2]$ converges to non-trivial limits if $\beta_N$ is chosen in a finer critical window around $\hat\beta_c=1$ \cite{CSZ19a}, more precisely,
\begin{equation} \label{eq:sigma}
	\beta_N^2 = \frac{1}{R_N} \bigg(1 + \frac{\theta + o(1)}{\log N}\bigg) \quad\quad
	\text{for a fixed $\theta\in \R$} \,.
\end{equation}
For the solution of the mollified 2D SHE, the same limiting second moments were computed earlier by Bertini and Cancrini \cite{BC98} in the critical window
\begin{equation}\label{eq:DBbeta2}
\beta^2_\epsilon= \frac{2\pi}{\log \frac{1}{\epsilon}} \Big(1+\frac{\theta + o(1)}{|\log \epsilon|}\Big).
\end{equation}
Subsequently, third moment \cite{CSZ19b} and higher moments \cite{GQT21} of $\cZ_{N;\, t}(\varphi)$ were all shown to converge to finite limits as $N\to\infty$. This strongly suggest that, with $\beta_N$ chosen in the critical window for some fixed $\theta\in \R$, we should have
\begin{equation} \label{eq:convergence?}
	\cZ^{\beta_N}_{N;\, t}(\varphi) = \int_{\R^2} \varphi(x) \, Z_{N, \beta_N}^{\omega}([Nt], [\sqrt{N} x]) \, \dd x
	\ \xrightarrow[\ N\to\infty \ ]{d} \ \int_{\R^2} \varphi(x) \, \mathscr{Z}_t(\dd x) \,
\end{equation}
for some limiting random measure $\mathscr{Z}_t(\dd x)$. The boundedness of the first moment $\bbE[\cZ^{\beta_N}_{N;\, t}(\varphi)]$ already implies the tightness of $(Z_{N, \beta_N}^{\omega}([Nt], [\sqrt{N} \cdot ]))_{N\in \N}$ as a sequence of random measures on $\R^2$. The key challenge is to show that there is a unique sub-sequential weak limit in \eqref{eq:convergence?}, which is our main achievement in \cite{CSZ23a} that constructs the {\em Critical 2D Stochastic Heat Flow}.

In order to state this result properly, we generalise the point-to-plane partition functions
from \eqref{eq:paf2} by considering \emph{point-to-point partition functions}
between time points $M \le N \in \N_0 = \{0,1,2,\ldots\}$ and space points $x, y\in\Z^2$:
\begin{equation} \label{eq:ZMN}
	Z^{\beta}_{M,N}(x,y) := \E \Big[
	e^{\sum_{n=M+1}^{N-1} \{\beta \omega(n,S_n) - \lambda(\beta)\}}
	\, \ind_{S_N = y} \,\Big|\, S_M = x \Big]  \,,
\end{equation}
and defining the process of rescaled random measures
$\cZ_N := \big(\cZ_{N;\, s,t}(\dd x, \dd y)\big)_{0\leq s\leq t}$ by
\begin{equation} \label{eq:Zappr}
 \cZ^{\beta_N}_{N;\, s,t}(\dd x, \dd y) :=
	N \, Z^{\beta_N}_{[ Ns], [ Nt ]}([\sqrt{N} x], [\sqrt{N}y]) \, \dd x \, \dd y \,.
\end{equation}
Point-to-point partition functions correspond to solutions of the SHE with general starting time
and delta initial conditions, while point-to-plane partition function corresponds to starting at time~$0$
from a constant~$1$ initial condition.

\begin{theorem}[Critical 2D SHF]\label{th:main0}
Let $\beta_N$ be chosen in the critical window \eqref{eq:sigma} for some $\theta\in \R$. As $N\to\infty$, the process of random measures
$\cZ^{\beta_N}_N = (\cZ^{\beta_N}_{N;\, s,t}(\dd x, \dd y))_{0\le s \le t}$ converges in finite dimensional distribution
to a {\it unique} limit $\mathscr{Z}^\theta = (\mathscr{Z}^\theta_{s,t}(\dd x, \dd y))_{0 \le s \le t < \infty}$
called the \emph{critical $2d$ Stochastic Heat Flow},
which does not depend on the law of the disorder $\omega$ (assuming $\bbE[\omega]=0$, $\bbE[\omega^2]=1$, with finite exponential moments).
The first and second moments are
\begin{align}
	\label{eq:mean-SHF}
	\bbE[\mathscr{Z}^\theta_{s,t}(\dd x , \dd y)]
	&= g_{t-s}(y-x) \, \dd x \, \dd y \,, \\
	\label{eq:var-SHF}
	\bbcov[\mathscr{Z}^\theta_{s,t}(\dd x ,\dd y),
	\mathscr{Z}^\theta_{s,t}(\dd x' , \dd y')]
	&= K_{t-s}^\theta(x,x'; y, y') \, \dd x
	\, \dd y \, \dd x' \, \dd y' \,,
\end{align}
where $g_t(\cdot)$ is the heat kernel on $\R^2$ and $K^\theta_{t-s}$ is defined by
\begin{equation}
\label{eq:m2-lim}
\begin{split}
	K_{t}^{\theta}(x,x'; y,y')
	&\,:=\, 4\pi \: g_{\frac{t}{2}}\big(\tfrac{y+y'}{2} - \tfrac{x+x'}{2}\big)
	\!\!\! \iint\limits_{0<a<b<t} \!\!\! g_{2a}(x'-x) \,
	G_\theta(b-a) \, g_{2(t-b)}(y'-y) \, \dd a \, \dd b \,,
\end{split}
\end{equation}
where $G_\theta(t) := \int_0^\infty \frac{e^{(\theta-\gamma)s} \, s \, t^{s-1}}{\Gamma(s+1)} \,\dd s$
is a so-called \emph{Volterra function}, linked to the \emph{Dickman subordinator}
\cite{CSZ19a}.
\end{theorem}

\begin{remark}
For simplicity, we have stated Theorem~\ref{th:main0} for directed polymer partition functions based on an aperiodic random walk on $\Z^2$, so that the normalisation
of the SHF matches the one of the mollified SHE. When dealing with the simple random walk,
as in \cite{CSZ23a}, one needs to take periodicity issues into account, e.g.\ to define the integer
parts in \eqref{eq:Zappr}. We skip these minor details in this presentation.
\end{remark}

\subsection{More historical background}
We provide here more details on the motivation and the results that led to the SHF.

As we already mentioned, the second moment of the LHS of \eqref{eq:convergence?}
 was first investigated by Bertini and Cancrini in the late 1990's
\cite{BC98} in the setting of the mollified 2D SHE,
exploiting the connection with the delta-Bose gas \cite{AGHKH05}.
They showed that the variance has a finite
non-zero limit as $N\to\infty$, corresponding to the kernel \eqref{eq:m2-lim},
under the scaling \eqref{eq:DBbeta2}.
However, convergence of the variance to a non-zero limit does not rule out
the possibility of a trivial (deterministic) scaling limit $\mathscr{Z}_t(\dd x) \equiv \dd x$.

The study of 2D SHE was mostly silent for several years, since the study of moments of order
higher than two is considerably more involved. In \cite{CSZ19a}
we obtained deeper insight on the second moment,
in the setting of directed polymers, via a \emph{probabilistic
approach} that gives a renewal representation for the second moment in terms of the
so-called \emph{Dickman subordinator}. This, in turn, allowed us to prove \emph{convergence of the
third moment of the LHS of \eqref{eq:convergence?} to a finite (explicit) limit},
see \cite{CSZ19b}. This guaranteed that any subsequential weak
limit $\mathscr{Z}_t(\dd x)$ in \eqref{eq:convergence?} must have the same covariance kernel \eqref{eq:m2-lim},
and hence must be truly random.

\begin{wrapfigure}[16]{r}{0.7\textwidth}
    \centering     
\hspace{1cm}

    \begin{tikzpicture}[scale=0.4]
        \foreach \i in {6,7,...,8}{
            \draw[-, thick] (2*\i,-6) -- (2*\i,6);
        }
        \draw[-,thick] (2,-6) -- (2,6);
        \draw[-,thick] (0,-6) -- (0,6);
        \draw[-,thick] (6,-6) -- (6,6);
        \draw[-,thick] (8,-6) -- (8,6);
        \draw[-,thick] (20,-6) -- (20,6);

        \fill[black] (2,0) circle [radius=0.175];
        \fill[black] (6,0) circle [radius=0.175];
        \fill[black] (8,2.5) circle [radius=0.175];
        \fill[black] (12,2.5) circle [radius=0.175];
        \fill[black] (14,-3) circle [radius=0.175];
        \fill[black] (16,-2) circle [radius=0.175];

        \draw[thick] (0,4) to [out=30,in=120] (8,2.5);
        \draw[thick] (0,2) to [out=0,in=120] (2,0);
        \draw[thick] (0,-2) to [out=0,in=-120] (2,0);
        \draw[thick] (0,-4) to [out=-20, in=200] (14,-3);
        \draw[thick] (6,0) to [out=75,in=200] (8,2.5);
        \draw[thick] (6,0) to [out=-30,in=200] (13,-1);
        \draw[thick] (13.5,-0.9) to [out=20,in=180] (20,1);
        \draw[thick] (12,2.5) to [out=-60, in=120] (14,-3);
        \draw[thick] (12,2.5) to [out=45,in=180] (20,3);
        \draw[thick] (16,-2) to [out=30,in=180] (20,-1);
        \draw[thick] (16,-2) to [out=-30,in=175] (20,-3);

        \draw[-,thick,decorate,decoration={snake, amplitude=.4mm,segment length=2mm}] (2,0) -- (6,0);
        \draw[-,thick,decorate,decoration={snake, amplitude=.4mm,segment length=2mm}] (8,2.5) -- (12,2.5);
        \draw[-,thick,decorate,decoration={snake, amplitude=.4mm,segment length=2mm}] (14,-3) -- (16,-2);

        \node at (0.5,-6) {{$0$}};
        \node at (2.5,-6) {{$a_1$}};
        \node at (6.5,-6) {{$b_1$}};
        \node at (8.5,-6) {{$a_2$}};
        \node at (12.5,-6) {{$b_2$}};
        \node at (14.5,-6) {{$a_3$}};
        \node at (16.5,-6) {{$b_3$}};
        \node at (20.5,-6) {{$1$}};

        \node at (0.7,4.9) {{$x^1$}};
        \node at (0.7,2.8) {{$x^2$}};
        \node at (0.7,-0.8) {{$x^3$}};
        \node at (0.7,-3.5) {{$x^4$}};

        \node at (2.7,-0.7) {{$x_1$}};
        \node at (6.7,-0.8) {{$y_1$}};
        \node at (8.7,1.7) {{$x_2$}};
        \node at (12.7,2) {{$y_2$}};
        \node at (14.7,-2) {{$x_3$}};
        \node at (16.7,-1) {{$y_3$}};

        \node at (4,1) {\tiny {$\{2,3\}$}};
        \node at (10,3.3) {\tiny {$\{1,2\}$}};
        \node at (15,-4) {\tiny {$\{1,4\}$}};
    \end{tikzpicture}
    \caption{Collision diagram for the fourth moment of the SHF.}
    \label{fig:moments}
    \end{wrapfigure}
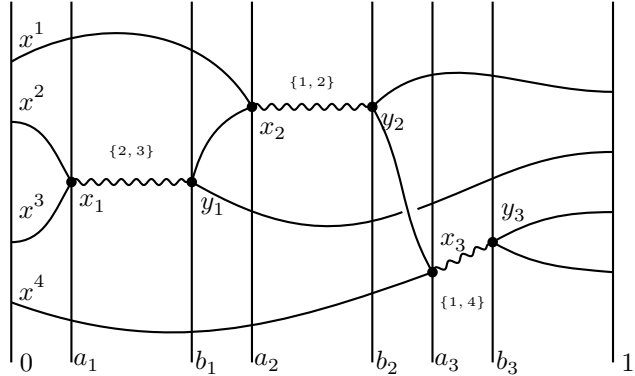

The key open problem left was then the uniqueness of subsequential limits, which
would yield the \emph{convergence to a unique limit in \eqref{eq:convergence?}}.
For this purpose, a natural idea is to compute all (integer) moments, beyond the first three. This was achieved
by Gu, Quastel and Tsai \cite{GQT21} by a Markovian operator approach
(assuming $L^2$ initial conditions). This led
to explicit \emph{formulas for all integer moments of the SHF}, which generalised the formulas
in \cite{CSZ19b}. These are represented as series of collision diagrams of independent 2D Brownian motions
with {\it point interactions} as in Figure~\ref{fig:moments},
where wiggle lines denote streams of pairwise collisions and are given weight $G_\theta(b-a) \, g_{\frac{b-a}{4}(y-x)}$,
while solid lines depict transition probabilities of the Brownian motions.
However, moments turn out to grow too fast to uniquely determine the distribution of the SHF.
We refer to the next section for more discussion and recent results
on moment asymptotics.

\setlength{\columnsep}{0.05pt}

Uniqueness of subsequential weak limits in \eqref{eq:convergence?}
was finally settled in \cite{CSZ23a}, where we proved Theorem~\ref{th:main0}.
The main difficulty of the proof was the lack of a characterisation of the limit,
which led us to develop a strategy based on \emph{coarse-graining} and \emph{Lindeberg principles},
which we describe in the next subsection. Moment estimates
played an important technical role in the application of the Lindeberg principle
in \cite{CSZ23a}: to this end, we exploited and refined
the operator approach from \cite{GQT21}.

Recently, Tsai gave an axiomatic characterisation of the SHF \cite{T24}, which he then applied to
obtain the SHF as the limit of solutions of the mollified 2D SHE.
We will discuss this in more detail in Section \ref{sec:properties}.

\subsection{Proof outline of Theorem~\ref{th:main0}} \label{S:sketch}
Let us sketch the proof of Theorem~\ref{th:main0} \cite{CSZ23a}. As anticipated,
the key difficulty is proving that a limit exists without a characterisation of the limit.
Our strategy is to show that, for fixed test functions $\varphi, \psi$,
the laws of $(\cZ^{\beta_N}_{N}(\varphi, \psi))_{N\in\N}$ form
a \emph{Cauchy sequence}, i.e.,
\begin{equation} \label{eq:close}
	\text{\emph{$\cZ_M^{\beta_M}(\varphi, \psi)$ and $\cZ_N^{\beta_N}(\varphi, \psi)$
	are close in distribution for large $M, N \in \N$}} \,.
\end{equation}
This was inspired by the work of Kozma \cite{Koz07} on the convergence of
loop erased random walk on a bounded domain in $\Z^3$, where a characterisation of the limit is still lacking.

To establish \eqref{eq:close}, the idea is to define
 for each $\eps\in (0,1)$ a \emph{coarse-grained partition function},
$\mathscr{Z}_{\epsilon}^{(\cg)}(\varphi, \psi | \Theta)$ with a similar
multi-linear structure as $\cZ^{\beta_N}_N(\varphi, \psi)$, see \eqref{eq:Zstart2},
which depends on a family of
\emph{coarse-grained random variables $\Theta$} that replace the
microscopic disorder variables $\eta$ in \eqref{eq:Zstart2}.
The coarse-grained partition function is meant to be a bridge between
partition functions of different sizes $\cZ_M^{\beta_M}$ and $\cZ_N^{\beta_N}$.  More precisely, we
can approximate
$\cZ_N^{\beta_N}(\varphi, \psi)$ via coarse-graining on time-space scale $(\eps N, \sqrt{\eps N})$
to show that
\begin{equation}\label{eq:ZcgL2}
\cZ^{\beta_N}_{N}(\varphi, \psi)
\,\overset{L^2}{\approx}
\, \mathscr{Z}_{\epsilon}^{(\cg)}(\varphi, \psi
| \Theta_{N, \eps})
\end{equation}
for a suitable family $\Theta_{N, \eps}$ of weakly dependent random variables,
depending on $N$ and~$\eps$, where the error in the $L^2$ approximation
is small in $L^2$ for small $\eps$, \emph{uniformly in large~$N$}.

We then prove \eqref{eq:close} by showing the following approximation:
\begin{equation}\label{eq:ZcgNM}
\mathscr{Z}_{\epsilon}^{(\cg)}(\varphi, \psi | \Theta_{N, \eps})
\,
\overset{\rm dist}{\approx}
\, \mathscr{Z}_{\epsilon}^{(\cg)}(\varphi, \psi | \Theta_{M, \eps})
\end{equation}
which is obtained by a suitable \emph{Lindeberg principle}, applied to
multilinear polynomials of weakly dependent random variables $\Theta$,
see \cite[Appendix A]{CSZ23a}.
The key technical ingredients to control the error in the Lindeberg principle
are uniform bounds on absolute moments $\bbE[|\Theta_{N, \eps}|^k]$ for $k > 2$,
where $\Theta_{N, \eps}$ turns out to be close to the original averaged
partition function (a sign of self-similarity).

\medskip
\section{Properties of the Critical 2D Stochastic Heat Flow.}
\label{sec:properties}

\begin{figure}[h]
\centering
\begin{minipage}[b]{.4\columnwidth}
\centering
\includegraphics[width=\columnwidth]{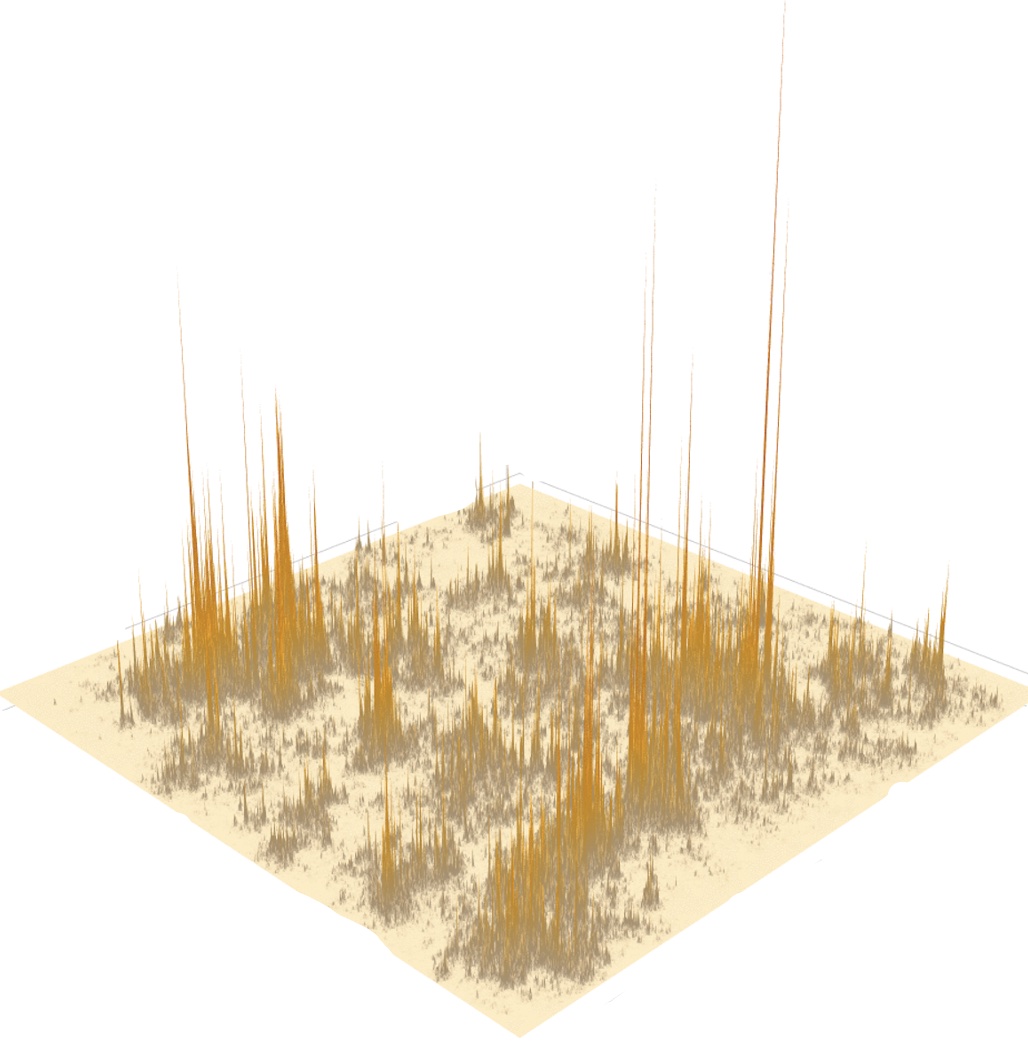}
\end{minipage}%
\begin{minipage}[b]{.03\columnwidth}
~
\end{minipage}%
\begin{minipage}[b]{.4\columnwidth}
\centering
\includegraphics[height=70mm]{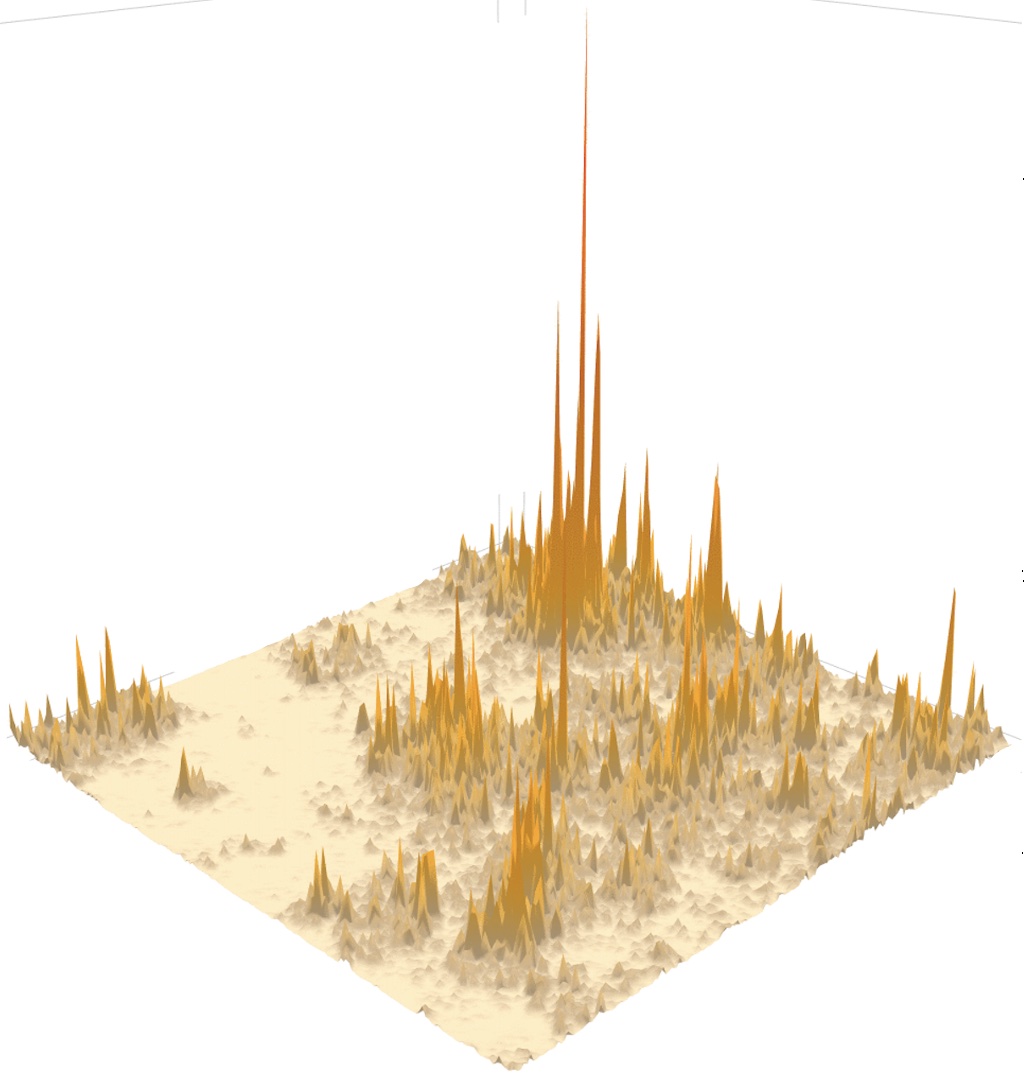}
\end{minipage}

\begin{minipage}[b]{.4\columnwidth}
\centering
\includegraphics[width=\columnwidth]{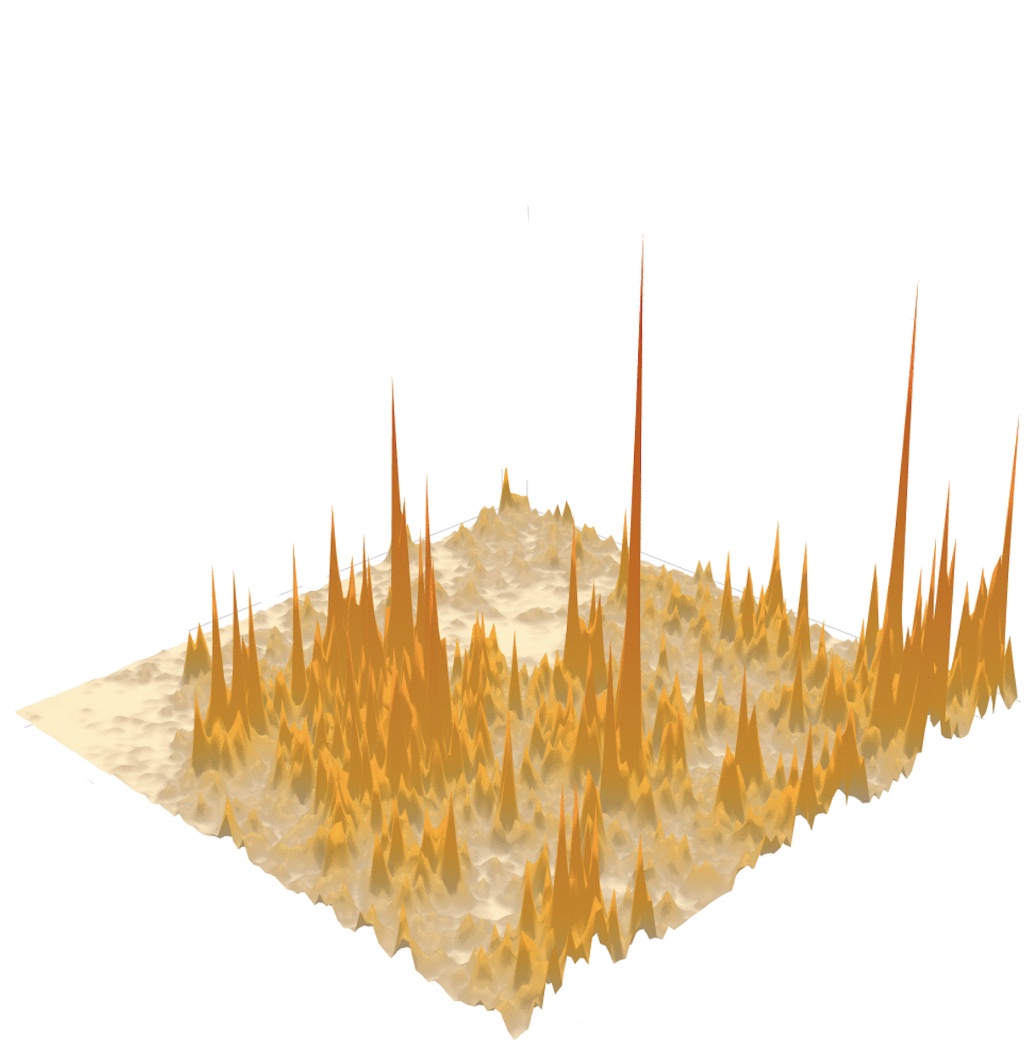}
\end{minipage}
\begin{minipage}[b]{.03\columnwidth}
~
\end{minipage}%
\begin{minipage}[b]{.4\columnwidth}
\centering
\includegraphics[width=\columnwidth]{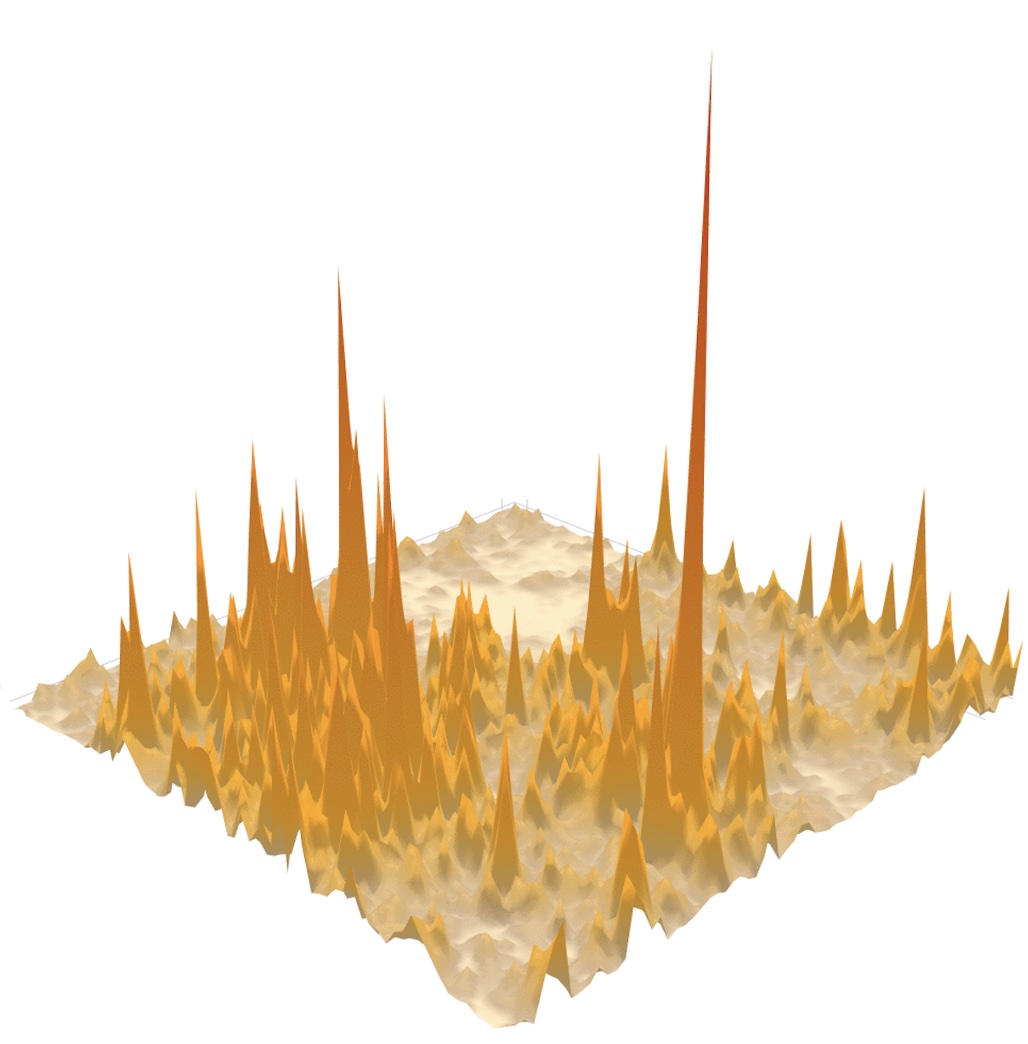}
\end{minipage}
\caption{
The top-left picture is a simulation of the Critical 2d SHF with $\theta=0$. The other three pictures
show how this looks when successively zoom into its scales, which effectively
corresponds to lowering the disorder parameter $\theta$, see the scaling
covariance property \eqref{eq:scaling}.} \label{image}
\end{figure}

Since the construction of the SHF, there have been many recent works exploring its properties.
We mention some of the highlights in this section. There is still much more to be explored.

\subsection{Axiomatic characterisation.}
An axiomatic characterisation of the Critical 2D SHF has been provided in \cite{T24} in the spirit of a Lindeberg principle.

\begin{theorem}[Characterisation of the Critical $2d$ SHF]\label{thm:tsai}
Let $Z=(Z_{s,t}(\cdot, \cdot))_{s\leq t}$ be a stochastic process taking values in
the space of locally finite non-negative measures $M_+(\R^2\times\R^2)$ on $\R^2\times\R^2$
equipped with the vague topology   such that:
\begin{itemize}
\item[(1)] for any $s<t<u$, the random measures $Z_{s,t}$, $Z_{t, u}$ and $Z_{s, u}$ satisfy a
form of Chapman-Kolmogorov property $Z_{s, u} = Z_{s, t} \, \bullet \, Z_{t,u}$;
\item[(2)] for any $s<t<u$, $Z_{s,t}$ and $Z_{t,u}$ are independent;
\item[(3)] for any $s<t<u$, $1\leq n\leq 4$, and $\phi_i,\psi_i\in L^2(\R^2)$ for $i=1,\ldots,4$, the mixed moments
$\bbE\big[\prod_{i=1}^n Z_{s,t}(\phi_i,\psi_i) \big]$ agree with that of the SHF with parameter $\theta\in\R$.
\end{itemize}
Then $Z$ has the same law as the Critical $2d$ Stochastic Heat Flow with parameter
$\theta$. Furthermore, $(Z_{s,t}(\cdot, \cdot))_{s\leq t}$ admits
a version that is almost surely continuous in $s\leq t$.
\end{theorem}
The Chapman-Kolmogorov type property (1)  in the above theorem for the SHF would formally write as
\begin{align*}
\mathscr{Z}^\theta_{s,u}(\mathrm{d} x, \mathrm{d} y) \, ``= " \int_{\R^2 }
\mathscr{Z}^\theta_{s,u}(\mathrm{d} x, \mathrm{d} z) \mathscr{Z}^\theta_{s,u}(\mathrm{d} z, \mathrm{d} y).
\end{align*}
However, the meaning of this expression is unclear since the SHF is a measure and so the integral
over $\mathrm{d} z$ is ill-defined.
A formalisation of the Chapman-Kolmogorov was provided in \cite{CM24}.

\subsection{On the martingale problem for the SHF.}
A martingale problem formulation of the Critical 2D SHF was given
in \cite{N25} and \cite{C25a}:
\begin{theorem}\label{thm: martingale}
    Let $\phi\in C^+_c(\bbR^2)$, $\psi\in C^2_b(\R^2)$ and $\theta\in\R$. Let
     $   \mathscr{Z}_t^{\theta,\phi}(\psi)
        :=\iint \phi(x) \, \mathscr{Z}_t^\theta(\mathrm{d} x, \mathrm{d} y) \, \psi(y) .$
    Then
    \begin{align}\label{eq:MG}
        \mathscr{M}_t^{\theta,\phi} (\psi)
        :=
        \mathscr{Z}_t^{\theta,\phi}(\psi)
        -\int \phi(x) \psi(y) \mathrm{d} x \, \mathrm{d} y
        -\int^t_0 \mathscr{Z}_s^{\theta,\phi} \big( \tfrac{1}{2}\Delta \psi \big) \mathrm{d} s
    \end{align}
    defines a continuous martingale starting at $0$ at time $0$ and with quadratic variation
    $\langle \mathscr{M}^{\theta,\phi}(\psi) \rangle_t$ obtained as a suitable limit via
    regularisation.
\end{theorem}

Although the martingales appearing in \eqref{eq:MG} are quite natural, identifying their quadratic variation is a non-trivial task and were only defined through a limiting procedure. The uniqueness of the martingale problem remains open, which probably requires a more explicit characterisation of the quadratic variation and specify its dependence on the disorder strength $\theta$.

\medskip

Moving from the foundational aspects of the SHF, we will now describe some of its qualitative properties.

\subsection{Translation invariance and scaling covariance.}
We have:
\begin{theorem}[Translation invariance and scaling covariance] \label{th:main1}
The Critical 2D SHF $(\mathscr{Z}_{s,t}^{\theta}(\mathrm{d} x , \mathrm{d} y))_{0 \le s \le t <\infty}$ is translation invariant in law:
\begin{equation*}
	(\mathscr{Z}_{s+\mathsf{a}, t+\mathsf{a}}^\theta(\mathrm{d} (x+\mathsf{b}) , \mathrm{d} (y+\mathsf{b})))_{0 \le s \le t <\infty}
	\stackrel{\rm dist}{=}
	(\mathscr{Z}_{s,t}^{\theta}(\mathrm{d} x , \mathrm{d} y))_{0 \le s \le t <\infty}
	\quad \ \forall \mathsf{a} \ge 0, \ \forall \mathsf{b} \in \R^2 \,,
\end{equation*}
and it satisfies the following scaling relation:
\begin{equation}\label{eq:scaling}
	(\mathscr{Z}_{\mathsf{a} s, \mathsf{a} t}^\theta(\mathrm{d} (\sqrt{\mathsf{a}} x) , \mathrm{d} (\sqrt{\mathsf{a}} y)))_{0 \le s \le t <\infty}
	\stackrel{\rm dist}{=}
	(\mathsf{a}\, \mathscr{Z}_{s,t}^{\theta+ \log \mathsf{a}}(\mathrm{d} x , \mathrm{d} y))_{0 \le s \le t <\infty}
	\quad \ \forall \mathsf{a} > 0 \,.
\end{equation}
\end{theorem}
The translation invariance is inherited from the
invariance of the white noise.
The scaling covariance property is illustrated in the simulation of the SHF in Figure~\ref{image}:
the top-left image is a picture of the SHF at $\theta=0$, while the following images are
successive zoom-ins and scaled up screenshots of the top-left image around its centre.
The zoom-ins appear to be smoother that the original picture due to the shift of the parameter $\theta$
by $\log \mathsf{a}$, which is negative for $\mathsf{a}\ll 1$ and thus shifts towards
the sub-critical regime.

\begin{figure}[h]
\centering
\begin{minipage}[b]{.52\columnwidth}
\centering
\includegraphics[width=\columnwidth]{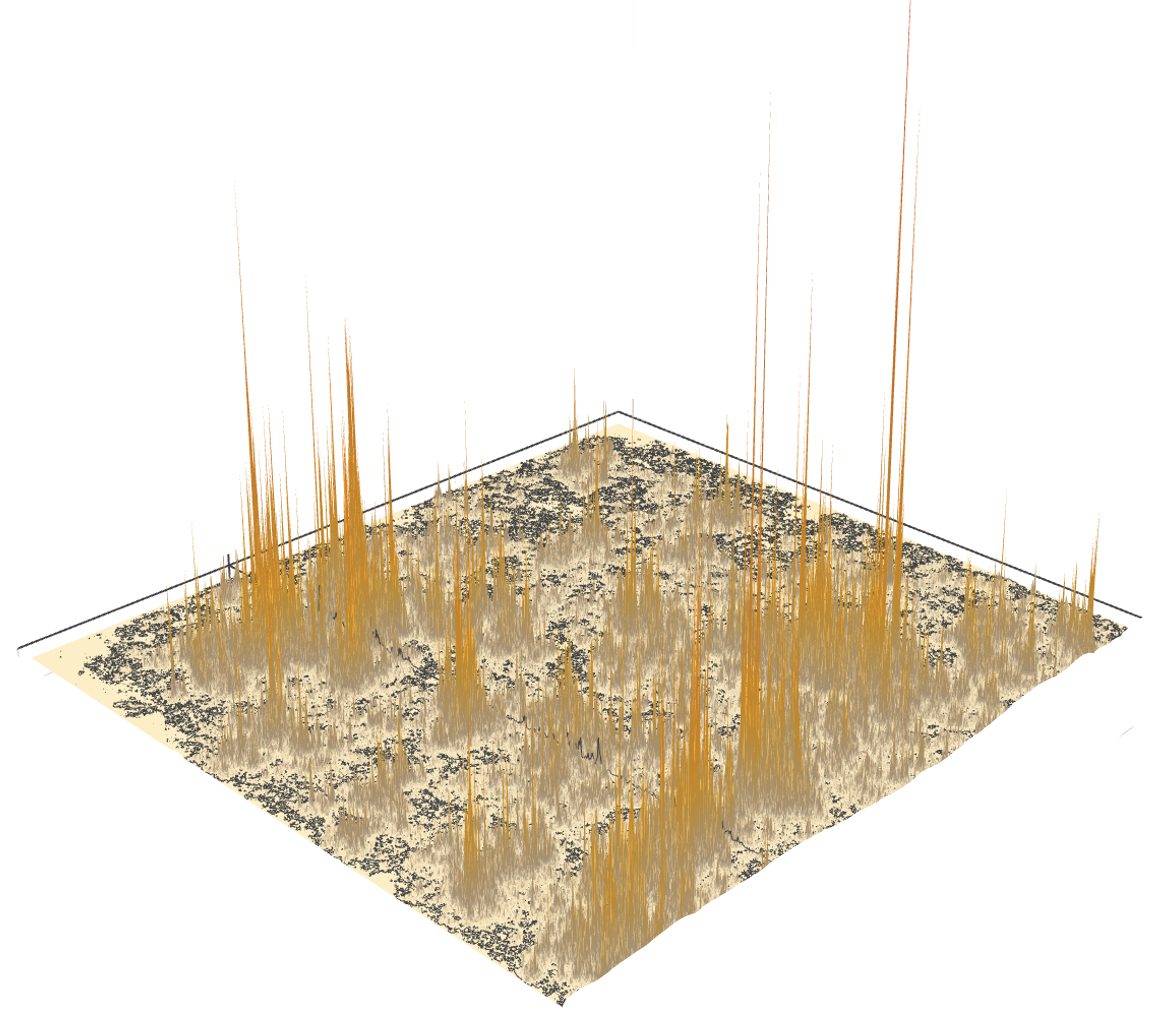}
\end{minipage}
\caption{The black contours show the level lines of the Critical 2D SHF corresponding to much lower values compared to the peaks.}\label{fig:level}
\end{figure}

\subsection{Singularity \& Regularity.}
The Critical 2D SHF is not a function:
\begin{theorem}[\cite{CSZ25}]\label{thm:sing}
For any $t>0$ and $\theta\in \R$, the marginal of the Critical 2D SHF  at time $t$,
$\mathscr{Z}^\theta_t(\mathrm{d}x):=\int_{y\in \R^2}\mathscr{Z}^\theta_t(\mathrm{d} y, \mathrm{d} x)$ is
a.s. singular with respect to the Lebesgue measure.
\end{theorem}

\vskip 2mm
However, it just fails to be a function:
\vskip 2mm
 \begin{theorem}[\cite{CSZ25}]\label{th:regularity-SHF}
Fix any $t > 0$ and $\theta\in\R$.
Almost surely, the critical $2d$ SHF
$\mathscr{Z}_{t}^\theta(\mathrm{d} x)$ belongs to $\mathcal{C}^{0-}
:= \bigcap_{\epsilon > 0} \mathcal{C}^{-\epsilon}$, where $\mathcal{C}^{-\epsilon}$ is the negative H\"older space
of order $-\eps$.
\end{theorem}

Nevertheless, since delta measures on $\R^d$ belong to $C^{-d}$ the above
 theorem indicates that it almost surely does not contains any atoms.
 The singularity of the SHF is apparent from its image in Figure~\ref{image}.
\vskip 2mm
The singularity of the SHF is obtained through another structural property, namely that
it is asymptotically log-normal when averaged over a small ball while the
 parameter $\theta$ is sent to $-\infty$
 in a precise way according to the radius of the ball.
Equivalently, we can send $\theta \to -\infty$ and adjust the ball radius accordingly,
in the following way.

\begin{theorem} \label{th:log-normality-SHF}
For a Euclidean ball $B(x,\delta) := \big\{y \in \R^2 \colon \ |y-x| < \delta \big\}$, denote the uniform distribution over it by
$\cU_{B(x,\delta)}(\cdot) := \frac{1}{\pi \delta^2} \, \mathrm{1}_{B(x,\delta)}(\cdot)$.
Then for any $t > 0$ and $x\in\R^2$, the following convergence in distribution holds:
\begin{equation}\label{eq:log-normal-SHF}
	\forall \rho \in (0,\infty) \colon \qquad
	\mathscr{Z}_{t}^{\theta} \big(\cU_{B(x,\sqrt{e^{\rho \theta}})}\big)
	\ \xrightarrow[\ \theta \to -\infty \ ]{d} \
	e^{\cN(0,\sigma^2) - \frac{1}{2}\sigma^2}
	\qquad \text{with } \sigma^2 = \log(1+\rho) \,.
\end{equation}
\end{theorem}

The regime $\theta \to -\infty$ means looking at the SHF in the weak disorder
limit. In this regime, a LLN and a CLT were obtained
in \cite{CCR25}, showing that the SHF displays Edwards-Wilkinson
fluctuations.

\begin{theorem} \label{th:EW-SHF}
As $\theta\to-\infty$ we have the convergence in distribution, for any $\varphi \in C_c(\R^2)$,
\begin{equation}\label{eq:main-SHF}
	\sqrt{|\theta|} \,
	\big\{ \mathscr{Z}_t^{\theta}(\varphi) - \bbE[\mathscr{Z}_t^{\theta}(\varphi)] \big\}
	\xrightarrow[\theta\to -\infty]{d} \, \mathcal{N} \big( 0 \,, v_{t,\varphi}\big)\,,
\end{equation}
where $v_{t,\varphi} :=
	\int_{\R^2 \times \R^2} \varphi(x) \, K_t(x,x') \, \varphi(x') \, \mathrm{d}x \, \mathrm{d}x'$
with
$K_t(x,x') := \int_{0}^{t} \frac{1}{2u} \, e^{-\frac{|x-x'|^2}{2u}} \, \mathrm{d}u$.
\end{theorem}

\subsection{Positivity of mass and local extinction.}
The picture of the Critical 2D SHF in Figure~\ref{image} might give the impression that there are
regions where the SHF puts no mass. However, this is not the case as shown in \cite{CT25,N25}:

\begin{theorem}\label{thm:mak}
For any initial condition $\phi\in C_c(\R^2)$, $\phi \ge 0$ 
which is not identically equal to $0$ we have that,
for any $t>0$, a.s.\ the SHF assigns strictly positive mass in any ball, i.e.,
$$
\mathscr{Z}^\theta_t(\phi, B(x,r))>0 \qquad
\text{ for all } r > 0, \ x\in \R^2 \,.
$$
\end{theorem}

The positivity of the mass of the SHF is depicted in Figure \ref{fig:level}. The back contours there indicate
the level lines corresponding to values which are much lower than the peaks. There appears to be a very
large separation of scales between the high and the low or even moderate values of the SHF
and it would be interesting to obtain a better understanding of the range of scales.

\medskip

On the other hand, the mass of the SHF converges locally to $0$ as time progresses as shown in \cite{CSZ25}:
\begin{theorem}
Let $\mathscr{Z}^\theta_t(A)$ be the mass of the Critical 2D SHF with (without loss of generality)
initial condition $1$, on a set $A\subset \R^2$. Then
\begin{equation} \label{eq:local-extinction}
	\text{for any bounded set $A\subset\R^2$, it holds that
	$\mathscr{Z}^\theta_t(A)\xrightarrow[t\to\infty]{d}{0}$}.
\end{equation}
\end{theorem}

\smallskip

The large-time regime $t \to \infty$ can be connected to the {\it strong disorder regime} $\theta\to+\infty$ via the scaling covariance property of the SHF \eqref{eq:scaling} which gives, in particular,
\begin{equation} \label{eq:scaling2}
	\mathscr{Z}_{t}^{\theta}(B(0,\sqrt{t}))
  \,\overset{d}{=}\,
 t \, \mathscr{Z}_{1}^{\theta + \log t}(B(0,1)) \,.
\end{equation}
It was recently proved in \cite{CT25} that $\mathscr{Z}^\theta_t(A) \to 0$ in distribution as
$\theta \to \infty$, which implies $t^{-1} \, \mathscr{Z}_{t}^{\theta}(B(0,\sqrt{t})) \to 0$ as
$t \to \infty$ by \eqref{eq:scaling2}. 
A quantitative strengthening was obtained in \cite{BCT25}, where optimal rates
were provided in both regimes $t\to\infty$ and $\theta \to \infty$ (possibly combined).

\begin{theorem}
There are constants $\delta > 0$ and $0 < c' < c'' < \infty$ such that the following  holds for any $t \ge 0$ and $\theta\in\R$:
  \begin{equation} \label{eq:SHF-largeball}
    \text{with probability at least $\,1 - \tfrac{1}{\delta} \, e^{-\delta \, t \, e^{\theta}}$:} \qquad
    \begin{cases}
      \mathscr{Z}_t^\theta\bigl(
    B\bigl(0,e^{c' \, t  \, e^{\theta}}
    \sqrt{t \,} \, \bigr)\bigr)
    \le t \, e^{-\delta \, t \, e^{\theta}} \,,
    \\
    \rule{0pt}{1.4em}
    \mathscr{Z}_t^\theta\bigl(
    B\bigl(0,e^{c'' \, t  \, e^{\theta}}
    \sqrt{t \,} \, \bigr))
    \ge t \, e^{\delta \, t \, e^{\theta}}  \,.
    \end{cases}
  \end{equation}
\end{theorem}

\subsection{Moment asymptotics and comparisons with GMC.}

Gaussian Multiplicative Chaos (GMC) is a distinguished, universal model in random geometries
with far reaching connections to Liouville Quantum Gravity,  random maps, turbulence and more. We refer to the reviews \cite{RV14,BP25,RV25} for comprehensive accounts.

Given a reference measure $\sigma$ on $\R^d$, the GMC is a random measure
whose density w.r.t.\ $\sigma$ is formally the 
exponential of a (typically generalised) Gaussian process 
$\big(\mathsf{X}(x) \big)_{x\in \R^d}$:
 \begin{align}\label{GMCdef}
 \mathsf{M}(\mathrm{d} x):= e^{ \mathsf{X}(x) -\frac{1}{2} \bbE[\mathsf{X}^2] } \,\sigma(\mathrm{d} x) \,.
 \end{align}
Some care is needed to make sense of this formula when $\mathsf{X}$ is distribution valued,
in which case $\mathsf{M}$ is typically \emph{singular} w.r.t.\ $\sigma$.
Certain assumptions are imposed on the covariance kernel $k(x,y)$ of $\mathsf{X}$
and the most studied one is that with a
logarithmic blow up $k(x,y)\approx \log \frac{1}{|x-y|}$ when $|x-y|\to0$.

It is a natural question whether the SHF $\mathscr{Z}^\theta(\dd x)$
coincides with a GMC on $\R^2$, call it $\mathsf{M}(\mathrm{d} x)$,
with reference measure $\sigma$ being the  Lebesgue measure,
for a suitable Gaussian field $\mathsf{X}$.
Equating  $
\E[ \mathsf{M}(\mathrm{d} x)  \mathsf{M} (\mathrm{d} y)] = \E[ \mathscr{Z}_t^\theta(\mathrm{d} x)  \mathscr{Z}_t^\theta (\mathrm{d} y)]$
determines the covariance of the candidate Gaussian field $\mathsf{X}$,
which must be
\emph{$\log\log$-correlated} (since the covariance of $\mathsf{M} = \mathscr{Z}_t^\theta$
 has logarithmic
blow-up). However, this cannot be the case as shown in \cite{CSZ23b} via comparison of the corresponding higher moments: for $g_\delta$ the heat kernel, the exists $\eta>0$ such that:
\begin{equation} \label{eq:asystrict}
	\liminf_{\delta\downarrow 0}  \;
	\frac{\bbE\big[ \mathscr{Z}_t^\theta(g_\delta)^h \big]}{\bbE\big[ \mathsf{M}(g_\delta)^h \big]}
	\ge 1+\eta > 1 \,, \qquad \text{for all $h\geq 3$}.
\end{equation}
Still, it is an open question how different is the SHF from a GMC.

In terms of moment comparisons, it follows by the definition of GMC
and the results in \cite{LiuZ24} that
\begin{align*}
\bbE\big[ \mathsf{M}(g_\delta)^h \big] = (C + o(1)) \, (\log \tfrac{1}{\delta})^{h\choose 2} \,,
\qquad
 \bbE\big[ \mathscr{Z}_t^\theta(g_\delta)^h \big]
 = (\log \tfrac{1}{\delta})^{{h\choose 2}+o(1)},\qquad \text{as $\delta\downarrow0$},
\end{align*}
which leaves open the question whether in small scales the SHF is a kind of perturbation of GMC.
On the other hand, if we keep $\delta$ fixed,
say $\delta = 1$, it was shown in \cite{GN25} that
\begin{align*}
	\bbE\big[ \mathscr{Z}_t^\theta(g_1)^h \big] \geq Ce^{e^h} \qquad \text{for all } h \in \N \,,
\end{align*}
which indicates that tail distribution of the SHF on a fixed ball is much heavier than that of GMC.

\subsection{SHF polymer measure and GMC on path space.}
The notion of GMC as in \eqref{GMCdef}
can be defined on more general measure spaces, such as the
path space $C\big([0,T];\R^d\big)$ with $\sigma$ being the Wiener measure. The exponential weight appearing in the Feynman-Kac
representation of the mollified SHE \eqref{eq:FK1} defines a GMC on such a path space,
and one may wonder whether this GMC structure is preserved in the limit of the SHF. 

It was proved in \cite{CM24} that the SHF, for any disorder parameter $\theta$,
can be lifted to a corresponding measure on the path space, but this measure
turns out not to be a GMC with respect to the Wiener measure (because the so-called
second moment measure is not absolutely continuous w.r.t.\ the product
of Wiener measure with itself).
Nevertheless, it was recently shown in \cite{CT25} that a {\it conditional GMC property}
still holds, namely the path measure corresponding to a parameter
strength $\theta$ can be represented as a GMC with respect to
the path measure with any smaller parameter $\theta'<\theta$, which
acts as a reference measure.
A predecessor of this phenomenon in polymer settings on hierarchical lattices
appeared in \cite{Cla23}.


\subsection{Black noise and noise sensitivity.}

The last feature we discuss is the \emph{black noise} property \cite{Tsi04a}
of the SHF, recently established in \cite{GT25},
which is a continuuum version of the noise sensitivity property (see below).
An interesting corollary is that
any coupling of the SHF $\mathscr{Z}^\theta$ and a white noise
$\xi$ (adapted to the same filtration) must be such that \emph{$\mathscr{Z}^\theta$
and $\xi$ are independent} (see \cite{HP24} for a related result for the
directed landscape). As a consequence, \emph{the SHF $\mathscr{Z}^\theta$
cannot be the solution of a SPDE driven by white noise}.

An enhanced form of \emph{noise sensitivity} was recently established in \cite{CD25},
which concerns the partition functions $\cZ_N$
that approximate the SHF $\mathscr{Z}^\theta$, see Theorem~\ref{th:main0}.
When we 
perturb disorder by independently resampling each variable $\omega(n,z)$
with a fixed (small) probability $\epsilon \in (0,1)$,
we obtain perturbed partition functions $\cZ_N^\epsilon$,
which are shown in \cite{CD25}
to become \emph{asymptotically independent of
$\cZ_N$} as $N\to\infty$.
A corollary is that the SHF $\mathscr{Z}^\theta$ arising as the limit of $\cZ_N$
is independent of the white noise $\xi$ that arises from the rescaled disorder.
This provides an alternative proof that the SHF $\mathscr{Z}^\theta$
cannot be the solution of a SPDE driven by~$\xi$.

\medskip



\end{document}